\documentclass[11pt]{amsart}
\usepackage{amscd,amssymb}
\usepackage{amsmath}
\usepackage{graphicx}
\usepackage{fancybox}
\usepackage{epic,eepic}
\usepackage{amstext}
\usepackage{pstricks}
\usepackage{pst-node,pst-plot,pst-coil}
\usepackage{amsthm}
\usepackage{amsmath}
\usepackage[utf8]{inputenc}
\usepackage[english]{babel}
\usepackage[square,numbers]{natbib}
\usepackage{booktabs}
\usepackage[hidelinks]{hyperref}
\usepackage{comment}
\usepackage{mathtools}
\usepackage{tikz,pgfplots}
\usepackage{tabularx}
\usepackage{subfloat}
\usepackage{caption}

\def\H{\mathcal{H}}
\def\V{\mathcal{V}}
\def\W{\mathcal{W}}
\def\U{\mathcal{U}}

\def\Ut{\widetilde{\mathcal{U}}}
\def\Vt{\widetilde{\mathcal{V}}}
\def\Xt{\widetilde{\Phi}}
\def\Xtl{\widetilde{\Phi}_\ell}
\def\Xtlm{\widetilde{\Phi}_\ell^m}
\def\Yt{\widetilde{\Psi}}
\def\Ytl{\widetilde{\Psi}_\ell}
\def\Ytlm{\widetilde{\Psi}_\ell^m}

\def\T{\mathcal{T}}
\def\l{\ell}

\def\Il{\mathcal{I}_\ell}
\def\Pl{\mathcal{P}_\ell}
\def\Pit{\widetilde{\Pi}}
\def\C{\mathcal{C}}
\def\Cl{\mathcal{C}_\ell}
\def\Clt{\mathcal{C}_{\ell,T}}
\def\Clm{\mathcal{C}_\ell^m}
\def\Clms{\mathcal{C}_\ell^{m,*}}
\def\Cltm{\mathcal{C}_{\ell,T}^m}
\def\Cltms{\mathcal{C}_{\ell,T}^{m,*}}
\def\Pit{\widetilde{\Pi}}
\def\Pit{\widetilde{\Pi}}

\def\linh{\mathrm{span}}
\def\supp{\mathrm{supp}}

\def\xtl{\tilde{\varphi}_\ell}
\def\xtk{\tilde{\varphi}_k}
\def\xtlm{\tilde{\varphi}_\ell^m}
\def\ytl{\tilde{\psi}_\ell}

\def\ytlm{\tilde{\psi}_\ell^m}
\def\bli{b_{\ell,j}}
\def\blk{b_{\ell,k}}
\def\blis{b_{\ell,j}^{*}}
\def\blim{b_{\ell,j}^m}
\def\blkm{b_{\ell,k}^m}
\def\blims{b_{\ell,j}^{m,*}}
\newcommand{\tnorm}[2]{{\left\lVert #1 \right\rVert}_{L^2(#2)}}
\newcommand{\tnormb}[2]{\Bigg\lVert #1 \Bigg\rVert_{L^2(#2)}}
\newcommand{\tnormf}[2]{{\| #1 \|}_{L^2(#2)}}

\newcommand{\vnorm}[2]{{\left\lVert #1 \right\rVert}_{\V(#2)}}
\newcommand{\vnormb}[2]{\Bigg\lVert #1 \Bigg\rVert_{\V(#2)}}
\newcommand{\vnormf}[2]{{\| #1 \|}_{\V(#2)}}
\newcommand{\tsp}[3]{{\left( #1\,,\,#2 \right)}_{L^2(#3)}}

\newcommand{\tspf}[3]{{( #1\,,\,#2 )}_{L^2(#3)}}

\newcommand\ds{\,\text{d}s}

\newcommand\dx{\,\text{d}x}

\newtheorem{theorem}{Theorem}[section]
\newtheorem{corollary}[theorem]{Corollary}
\newtheorem{lemma}[theorem]{Lemma}
\newtheorem{assumption}[theorem]{Assumption}

\theoremstyle{definition}

\theoremstyle{remark}
\newtheorem{remark}[theorem]{Remark}
\numberwithin{theorem}{section}
\numberwithin{equation}{section}
\numberwithin{table}{section}
\numberwithin{figure}{section}

\allowdisplaybreaks

\begin{document}
\title[]{Multi-resolution Localized Orthogonal Decomposition for Helmholtz problems}
\author[]{Moritz Hauck$^\dagger$, Daniel Peterseim$^\dagger$}
\address{${}^{\dagger}$ Institute of Mathematics, University of Augsburg, Universit\"atsstr.~12a, 86159 Augsburg, Germany}
\email{\{moritz.hauck, daniel.peterseim\}@uni-a.de}
\thanks{The work of the authors is part of a project that has received funding from the European Research Council (ERC) under the European Union's Horizon 2020 research and innovation programme (Grant agreement No.~865751).}

\maketitle
\begin{abstract}
	We introduce a novel multi-resolution Localized Orthogonal Decomposition (LOD) for time-harmonic acoustic scattering problems that can be modeled by the Helmholtz equation. The method merges the concepts of LOD and operator-adapted wavelets (gamblets) and proves its applicability for a class of  complex-valued, non-hermitian and indefinite problems. It computes hierarchical bases that block-diagonalize the Helmholtz operator and thereby decouples the discretization scales. Sparsity is preserved by a novel localization strategy that improves stability properties even in the elliptic case. 
	We present a rigorous stability and a-priori error analysis of the proposed method for homogeneous media. In addition, we investigate the fast solvability of the blocks by a standard iterative method. 
	A sequence of numerical experiments illustrates the sharpness of the theoretical findings and demonstrates the applicability to scattering problems in heterogeneous media.
\end{abstract}

{\tiny {\bf Keywords.}~multi-resolution decomposition, multi-level method, numerical homogenization, Helmholtz problem,\linebreak 
	 \hphantom{\indent {\tiny {\bf Keywords.}}~} wave propagation, heterogeneous media
}\\
\indent
{\tiny {\bf AMS subject classifications.}  
	{\bf 65N12}, 
	{\bf 65N15},
	{\bf 65N30}, 
	{\bf 35J05}
}

\section{Introduction}


The concept of Localized Orthogonal Decomposition (LOD) has first been introduced in \cite{Malqvist2014,Henning2013} for elliptic problems. 
It is an approximately orthogonal (w.r.t. the energy inner product) decomposition of the energy space into a finite-dimensional, problem-adapted, mesh-based space with locally supported basis functions and an infinite-\linebreak dimensional remainder space. Using the finite-dimensional space as ansatz space for a Galerkin method (and allowing a moderate increase of the support of the basis functions), yields an optimally convergent multi-scale method that is capable of approximating solutions corresponding to arbitrarily rough $L^\infty$-coefficients, without any preasymptotic effects. 


This method has also proved to posses stabilizing effects for high frequency wave scattering applications \cite{Gallistl2015,Peterseim2016a,Peterseim2017a,Li2018}. For Helmholtz problems, the LOD method yields faithful approximations under the minimal resolution condition that the mesh size is coupled linearly to the wave number $\kappa$, and, if additionally, the support of the basis functions is allowed to increase logarithmically with $\kappa$. 
Its intrinsic ability to deal with high heterogeneity makes the LOD very appealing for such Helmholtz problems; see \cite{Brown2017,Verfuerth2020,Maier2020}. 


Recently, the numerical homogenization approach has been extended to multi-resolution approximation schemes that allow for a  decoupling into a hierarchy of discrete scales  \cite{Owhadi2017a}.  This concept has been popularized under the name \emph{gamblets}, and, up to now, has mostly been studied for the  homogenization of elliptic operators \cite{Owhadi2017,Owhadi2019}.



By bridging the gamblet approach with the well-developed LOD framework, we derive a multi-resolution (multi-level) LOD method that is applicable to a wide class of (complex-valued) non-hermitian and indefinite problems. As problem specific parameters make an unified (parameter-explicit) analysis nearly impossible, we  henceforth consider the Helmholtz problem for demonstration purposes.
Computing two sets of hierarchical bases for trial and test spaces, respectively, allows to block-diagonalize the operator. Note that, in contrast to the elliptic case, one hierarchical basis for both, trial and test space, does not suffice. 



This paper features a $\kappa$-explicit  stability and a-priori error analysis of the proposed multi-resolution method. 
For the elliptic case, gamblets suffer from severe instability, whenever the mesh size of the coarsest level is decreased; see also \cite{Maier2020a}. 
This is especially problematic for Helmholtz problems, where the mesh size of the coarsest level is coupled to $\kappa$ (for numerical illustrations, see Section \ref{subsec:improstab}). Using a novel basis construction enables us to cure this stability issue. The novel construction yields significant improvements also for single-level LOD methods as in \cite{Maier2020a,Maier2020b} for elliptic problems. We prove optimal rates of convergence of the proposed method, if a moderate increase of the support of the basis functions is allowed.   


The proposed multi-resolution method condenses the negative features of the linear system of equations arising from the discretization of the Helmholtz problem (indefinite,  deteriorating condition for increasing $\kappa$) to the first (small) block. For the remaining blocks, we prove that they are coercive and can be solved within a fixed number of \linebreak GMRES iterations, independent of levels, mesh sizes, and $\kappa$. This also inspires multi-level solver techniques for the Helmholtz equation, similarly as in \cite{Owhadi2019} for elliptic problems.



For the sake of simplicity, the analysis covers only Helmholtz problems in homogeneous media. Nevertheless, the theoretical results  carry over also to the case of heterogeneous media (for numerical experiments, see Section \ref{subsec:varcoeff}). However, the analytical well-posedness of the heterogeneous Helmholtz problem is problematic. In recent years, a lot progress has been made in this field; see  \cite{Sauter2018,Graham2018,Graham2019}. We believe that the results of the paper are also relevant for the simulation of other time-harmonic wave phenomena such as \cite{Brown2016,Gallistl2018,Maier2020}.

The remaining part of the paper is outlined as follows. Section \ref{sec:Helmhprop} defines the Helmholtz problem and recalls some of its fundamental properties. Section \ref{sec:notation} introduces the hierarchy of meshes, the Haar basis and corresponding mesh-based operators that will be the basis for the derivation of a prototypical multi-resolution method in Section \ref{sec:idealmethod}. Sections \ref{sec:localizetion} and \ref{sec:pracmethod} will then turn this ideal approach into a feasible practical method including a rigorous stability and a-priori error analysis. In Section \ref{sec:fastsolvers}, we prove that all blocks besides the first one can be inverted easily. Finally, Section \ref{sec:numexp} illustrates the performance of the method in a sequence of numerical experiments.

\section{Model Helmholtz problem}\label{sec:Helmhprop}
We consider the Helmholtz problem on a bounded Lipschitz domain $D\subset \mathbb{R}^d$, ($d = 1,2,3)$, which is scaled to unit size
\begin{subequations}\label{eq:Helmholtz}
\begin{equation}
-\Delta u  - \kappa^2 u = f\qquad \text{ in } D
\end{equation}
with boundary conditions of Dirichlet, Neumann, and Robin type
\begin{align}
u &= 0\qquad \text{ on } \Gamma_\mathrm{D},\\
\nabla u \cdot \nu &= 0\qquad \text{ on } \Gamma_\mathrm{N},\\
\nabla u \cdot \nu - i\kappa u &= 0\qquad \text{ on } \Gamma_\mathrm{R}.
\end{align}
\end{subequations}
Here, $\kappa\in(\kappa_0,\infty)$ with  $\kappa_0>0$ denotes the wave number, $i$ is the imaginary unit and $\nu$ is the outer unit normal vector. The right-hand side $f$ is assumed to be in $L^2(D)$ (the space of complex-valued square-integrable functions on $D$). It is also possible to use other types of boundary conditions like perfectly matched layers; see \cite{Chaumont-Frelet2018}. We assume a decomposition of the boundary into closed components
\begin{equation*}
\partial D = \Gamma_\mathrm{D} \cup \Gamma_\mathrm{N} \cup \Gamma_\mathrm{R},
\end{equation*}
where the intersection of the interior of the components is supposed to be pairwise disjoint. We allow $\Gamma_\mathrm{D}$ and $\Gamma_\mathrm{N}$ to be empty, however, for $\Gamma_\mathrm{R}$, we suppose a positive surface measure, i.e., 
\begin{equation}\label{eq:possurfmeasure}
	|\Gamma_\mathrm{R}|_{d-1}>0.
\end{equation}
Given the Sobolev space $W^{1,2}(D)$ (the space of complex-valued square-integrable functions on $D$ with square-integrable weak derivatives), we define the subspace
\begin{equation*}
\V := \left\{v \in W^{1,2}(D):\,v\vert_{\Gamma_\mathrm{D}} = 0\right\},
\end{equation*}
which is endowed with the usual  $\kappa$-dependent norm
\begin{equation*}
\vnorm{u}{D} := \sqrt{\tnorm{\nabla u}{D}^2 + \kappa^2 \tnorm{u}{D}^2}.
\end{equation*}
The weak formulation of \eqref{eq:Helmholtz} seeks $u\in\V$ such that for all $v\in \V$
\begin{equation}\label{eq:wfcont}
a(u,v) = \tsp{f}{v}{D}
\end{equation}
with the sesquilinear form $a:\V\times\V\rightarrow \mathbb{C}$ 
\begin{align*}
a(u,v) := \tsp{\nabla u}{\nabla v}{D} - \kappa^2\tsp{u}{v}{D} - i\kappa\tsp{u}{v}{\Gamma_\mathrm{R}}.
\end{align*}
Here, $\tsp{u}{v}{D}:=\int_D u \cdot \overline{v}\dx$ ($\,\overline{\,\cdot\,}$ is the complex conjugation) denotes the inner product of scalar or vector-valued functions in $L^2(D)$.  Similarly, $\tsp{u}{v}{\Gamma_\mathrm{R}}:= \int_{\Gamma_\mathrm{R}} u\overline{v}\ds$ (integration with respect to the surface measure) is the inner product of the space $L^2(\Gamma_\mathrm{R})$. The sesquilinear form $a$ is continuous, i.e., 
\begin{equation}\label{eq:conta}
|a(u,v)|\leq c\, \vnorm{u}{D}\vnorm{v}{D}
\end{equation}
with a generic constant $c>0$ independent of $\kappa$. In the following, we denote the real and the imaginary parts of a complex number by $\mathfrak{R}$ and $\mathfrak{I}$, respectively. 

The well-posedness of weak formulation \eqref{eq:wfcont} depends on the shape of $D$ and the choice of boundary conditions. The presence of Robin boundary conditions (c.f. \eqref{eq:possurfmeasure}) ensures unique solvability, i.e., for all $f \in L^2(D)$, the solution $u\in \V$ satisfies
\begin{equation}\label{polystab}
\vnorm{u}{D} \leq c_{\mathrm{stab}}(\kappa) \tnorm{f}{D},
\end{equation}
For the stability and a-priori error analysis in Section \ref{sec:pracmethod}, we assume a polynomial $\kappa$-dependence of $c_{\mathrm{stab}}$ (c.f. Assumption \ref{assumption:poly}), i.e., there exists a constant $c>0$ and $n\in \mathbb{N}$ such that
\begin{equation}\label{eq:poly}
c_{\mathrm{stab}}\leq c \kappa^n.
\end{equation} 
Clearly, condition \eqref{eq:poly} does not hold in general; see counter-examples with an at least exponential-in-$\kappa$ growth of $c_{\mathrm{stab}}$ for trapping domains \cite{Betcke2011}.

Nevertheless, \eqref{eq:poly} can be proven under certain geometric assumptions. In \cite{Melenk1995}, condition \eqref{eq:poly} is proved with $n = 0$ for bounded star-shaped domains with smooth boundary or bounded convex domains. Among the known admissible setups is also the case of pure Robin
boundary conditions ($\Gamma_\mathrm{R} = \partial D$) on Lipschitz domains \cite{Esterhazy2011}. Another possible example are truncated exterior Dirichlet problems, where the Sommerfeld radiation condition is approximated by truncating the (unbounded) exterior domain and applying Robin boundary conditions on the artificial boundary \cite{Hetmaniuk2007}.

\begin{remark}[Variable coefficients]
	In \cite{Sauter2018,Graham2018,Graham2019} the well-posedness of the Helmholtz problem in heterogeneous media is analyzed; see \eqref{eq:varcoeff} for the heterogeneous problem.
\end{remark}

An immediate consequence of \eqref{polystab} is the inf-sup stability of $a$
\begin{equation*}
 0<\alpha^{(1)}\leq \inf_{u\in \V} \sup_{v\in\V} \frac{\mathfrak{R} a(u,v)}{\vnorm{u}{D}\vnorm{v}{D}} = \inf_{v\in \V} \sup_{u\in\V} \frac{\mathfrak{R} a(u,v)}{\vnorm{u}{D}\vnorm{v}{D}}
\end{equation*}
with $\alpha^{(1)} = (2c_{\mathrm{stab}}(\kappa)\kappa)^{-1}$; see \cite{Melenk1995}. 

\section{Preliminaries}\label{sec:notation}
{\color{black}
This section lays the foundations for the construction of the multi-resolution LOD in the following sections. 
}


\subsection{Mesh hierarchy and Haar basis}

{\color{black} In what follows, we define a Haar basis given a sequence of nested meshes. The Haar construction is an essential ingredient for the construction of the proposed method as it induces its multi-resolution structure. Using the orthogonality properties of the Haar basis, a decoupling of discretization scales can be  achieved. }

Let $\T_1$ denote a (coarse) Cartesian mesh of $D$ with mesh size $H_1$ and let $\{\T_\l\}_{\l=1,...,L}$, $L\in \mathbb{N}$ be a hierarchy of meshes obtained by successive uniform refinements of $\T_1$. The mesh size at level $\l$ is given by $H_\l:=2^{-\l+1}H_1$. 
The patch of a union of elements $S$ in $\T_\l$ is defined as
\begin{equation*}
N^m_\l(S):= N_\l(N_\l^{m-1}(S)),\quad N_\l^0(S) := S, \quad N_\l(S) := \bigcup_{\T_\l\ni T:\, T\cap S\neq \emptyset} T.
\end{equation*} 

In the setting of a Cartesian mesh hierarchy, the corresponding Haar wavelet basis can be given explicitly. Denoting with $\mathbf{1}_S$ the characteristic function of the set $S$, we define
\begin{equation*}
\chi^{(0)}(x) := \mathbf{1}_{[0,1]}(x)\quad \text{ and }\quad\chi^{(1)}(x):= \chi^{(0)}(2x)-\chi^{(0)}(2x-1).
\end{equation*} 
The Haar basis can be written as $\mathcal{H}:=\bigcup_{\ell = 1,...,L} \mathcal{H}_\l$ with  $\mathcal{H}_1:=\{|T|^{-1/2}\mathbf{1}_T, T\in \T_1\}$ and 
 \begin{equation*}
\mathcal{H}_\l:=\left\{|T|^{-1/2}\chi^{(j_1)}(H_{\l-1}^{-1} x_1-a_{T,1})\cdot...\cdot\chi^{(j_d)}(H_{\l-1}^{-1} x_d-a_{T,d}) \,\Big|\, T\in \T_{\l-1},j \in \{0,1\}^d\backslash \mathbf 0\right\}
\end{equation*}
for $\l\geq 2$. Here, $a_{T,k}$ denotes the $k$-th coordinate of the bottom, front, left corner of $T\in \T_\l$ and $\mathbf 0$ is the $d$-dimensional zero vector. It is straightforward that  $\bigcup_{\l = 1,...,L} \H_\l$ is an othonormal basis of $\mathbb{Q}_0(\T_L)$. Henceforth, let $\phi_{\l,j}$, $j\in \{1,...,N_\l\}$ with $N_\l:=\#\H_\l$ denote the elements of $\H_\l$.
\begin{figure}[h]
	\begin{tabularx}{\textwidth}{XXXX}
	\includegraphics[width=\linewidth]{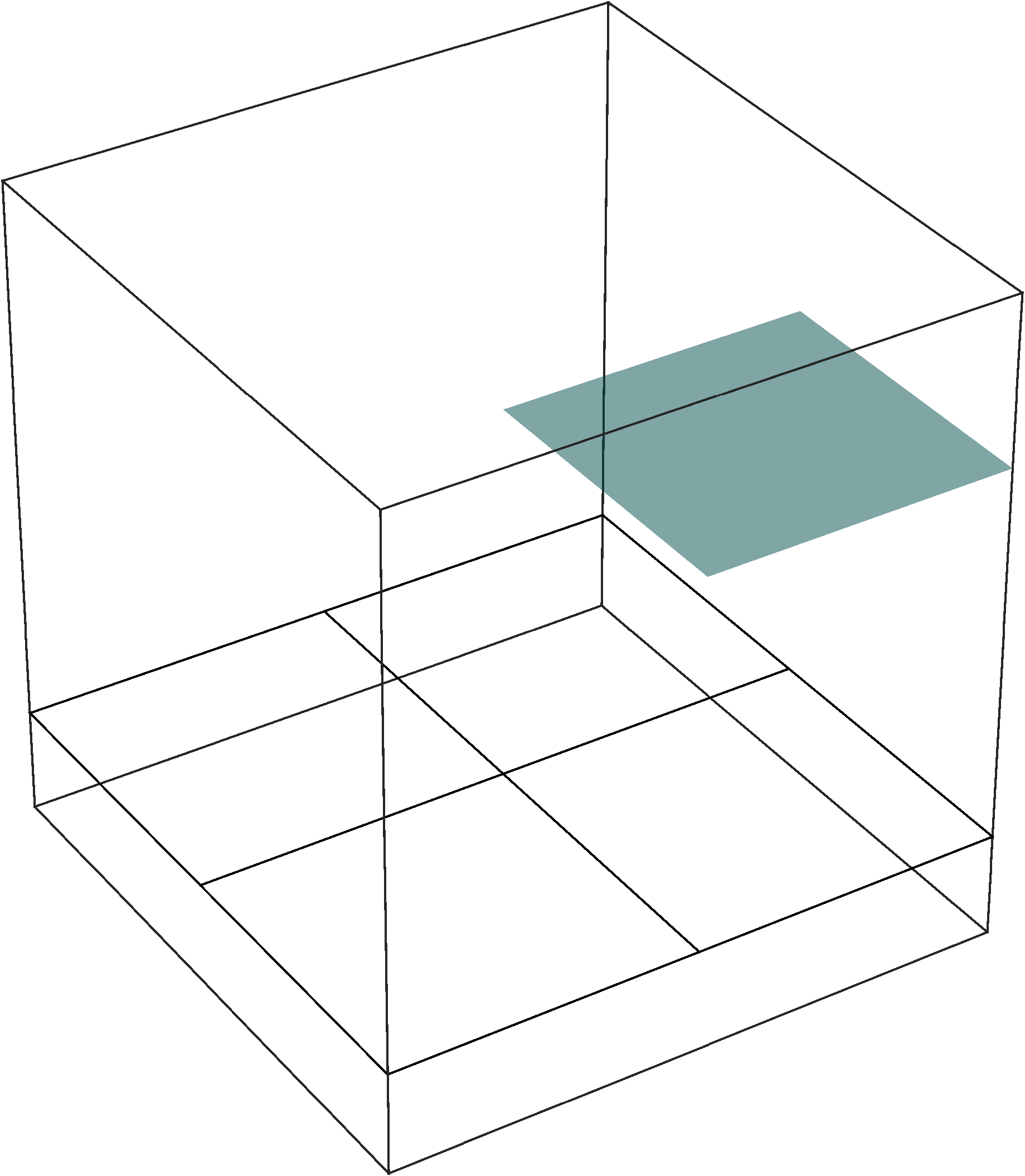}&
	\includegraphics[width=\linewidth]{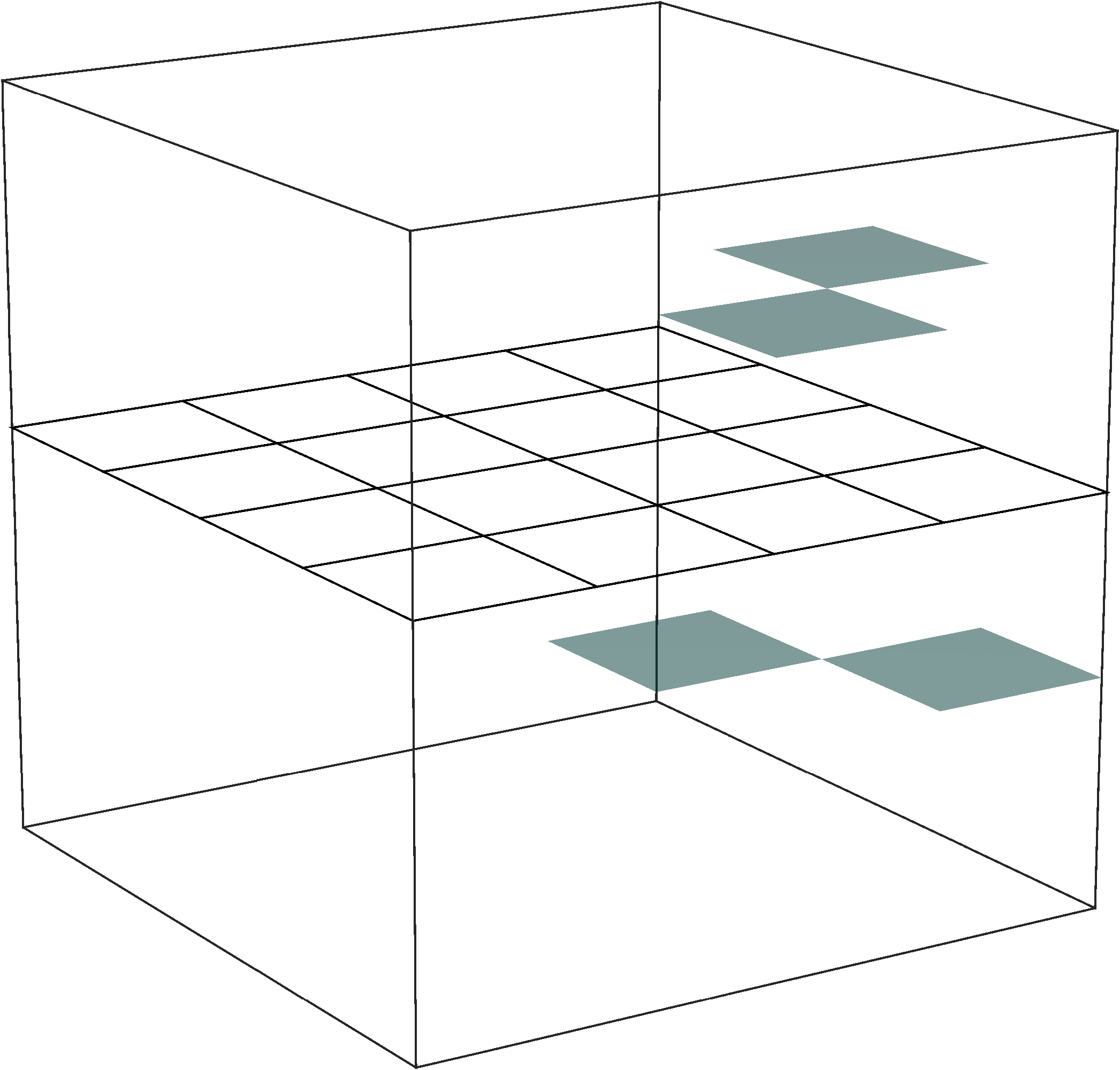}&
	\includegraphics[width=\linewidth]{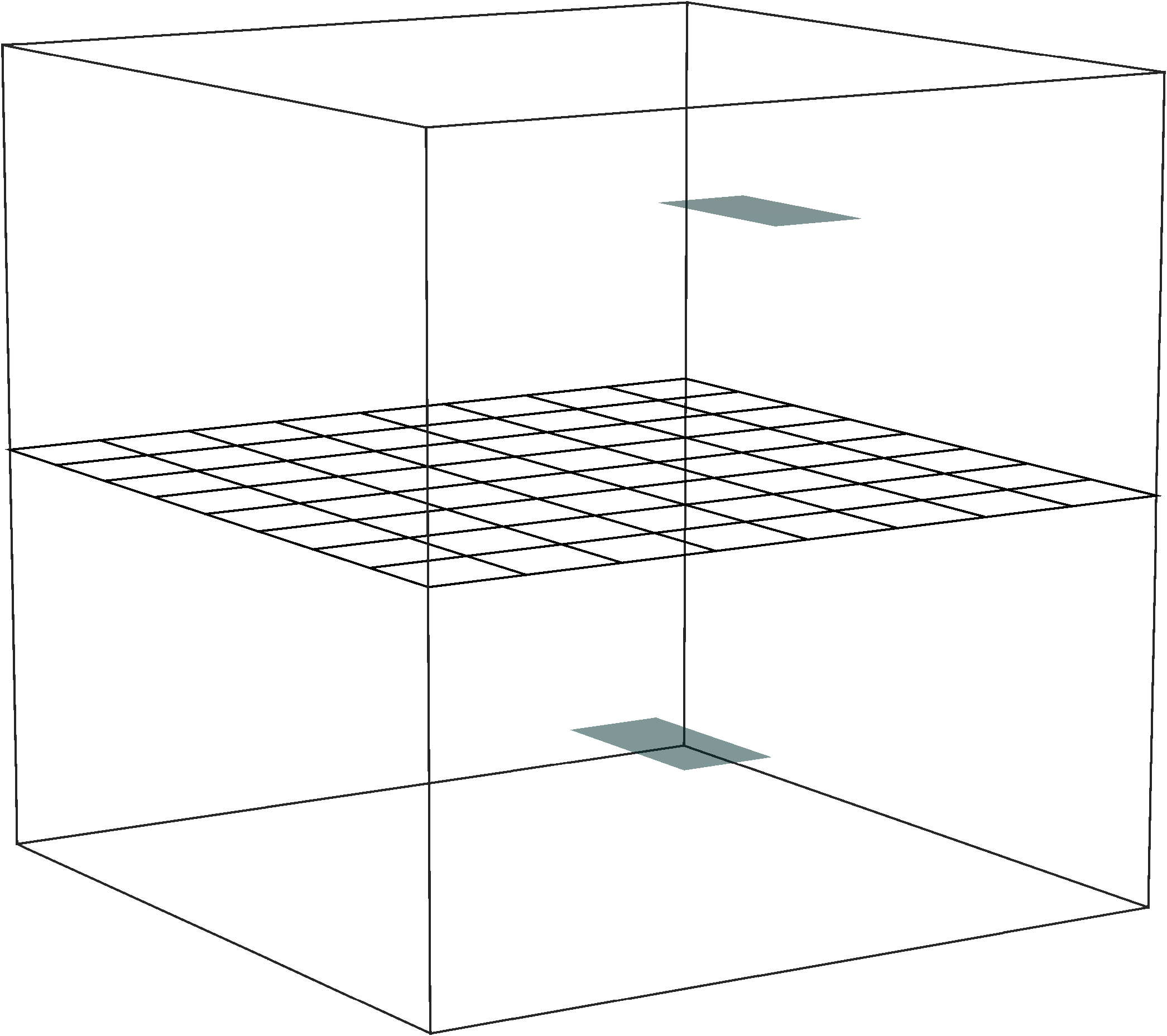}&
	\includegraphics[width=\linewidth]{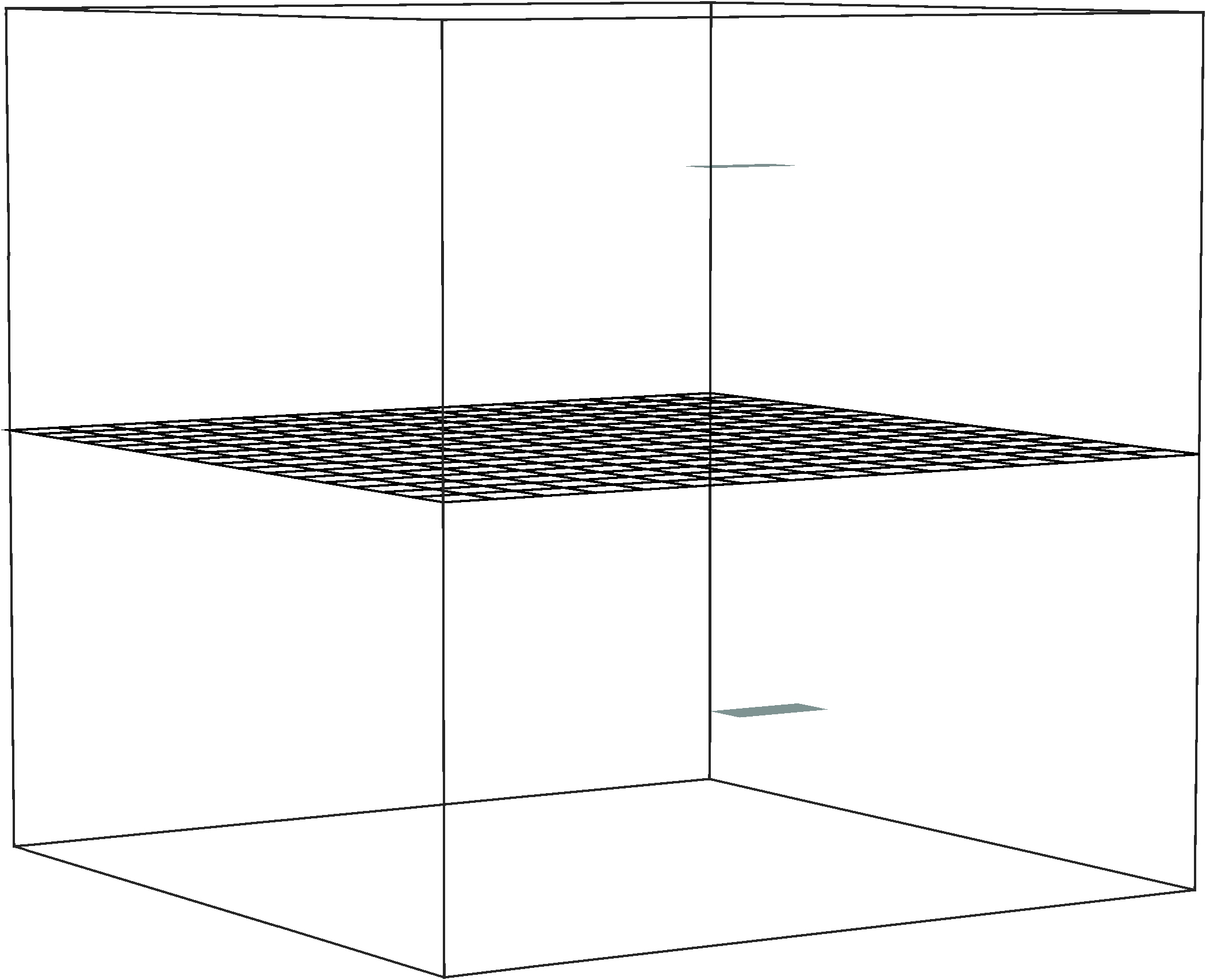}
	\end{tabularx}
\caption{Let $\T_1$ be a Cartesian mesh of the unit square with mesh size $H_1 = 1/2$. The meshes $\T_\l$ for $\l = 2,3,4$ are obtained by successive uniform refinements of $\T_1$. The pictures illustrate one (randomly chosen) Haar basis function per level.}
\label{fig:haarbasis}
\end{figure}

\subsection{Bubble functions} 
Clearly, the Haar spaces are non-conforming, i.e., $\H_\l\not\subset \V$. For the construction of the problem-adapted hierarchical basis, we introduce conforming companions to element-wise constant functions with the same element-averages.

For $\l = 1,...,L$ and all $T\in \T_\l$ define
\begin{equation*}
b_T := b(H_\l^{-1} x_1-a_{T,1})\cdot...\cdot b(H_\l^{-1} x_d-a_{T,d})\quad \text{ with }\quad b(x):= 6x(1-x)\mathbf{1}_{[0,1]}(x)
\end{equation*} 
and $a_{T,k}$ as above (note that $\int_{0}^{1}b(x)\dx = 1$). We define the linear mapping $\widetilde\Pi_\ell: \mathbb{Q}_0(\mathcal{T}_\ell)\rightarrow \mathcal{V}$ uniquely by setting
\begin{equation*}
\Pit_\l\mathbf{1}_T:=b_T\quad \text{for}\quad T\in\T_\l.
\end{equation*}
It maps $\mathbb{Q}_0(\T_\l)$-functions to conforming bubble functions with the same element averages. The operator is the right inverse of $\Pi_\ell$, where $\Pi_\l:L^2(D)\rightarrow\mathbb{Q}_0(\T_\l)$ denotes the $L^2$-orthogonal projection onto the space of element-wise constants with respect to $\T_\l$. The operators $\Pi_\l$ and $\Pit_\l$ satisfy for all $T\in\T_\l$
{\color{black}
	\begin{alignat}{2}
	\tnorm{\Pi_\l v}{T}&\leq \tnorm{v}{T},\quad\, \tnorm{(1-\Pi_\l)v}{T}&&\leq \pi^{-1}H_\l \tnorm{\nabla v}{T},\label{eq:Pi}\\
	\|\widetilde \Pi_\ell v\|_{L^2(T)}&\leq \|v\|_{L^2(T)},\qquad\quad \|\nabla \widetilde \Pi_\l v\|_{L^2(T)}&&\leq c_\mathrm{inv} H_\ell^{-1}\|v\|_{L^2(T)},\label{eq:Pit}
	\end{alignat}
	where $c_\mathrm{inv}>0$ is the constant from the inverse inequality (see e.g. \cite[Lemma 1.44]{DiPietro2012}).}

\begin{remark}[Construction of Haar basis and bubbles]
	The construction of the Haar basis and the bubbles is not restricted to Cartesian meshes. It can be performed for non-structured quadrilateral/hexahedral or simplicial meshes; see \cite{Alpert1993,Feischl2020}. If one considers tensor-product elements, higher order versions can be constructed; see \cite{Alpert1993,Maier2020a}.
\end{remark}
\begin{remark}[Tilde notation]
	Throughout the article, we write $a\lesssim b$ in short for $a \leq  cb$ with a constant $c$ independent of $\kappa$, $\l$, $H_\l$, $L$, and $m$ ($\,\gtrsim$ is defined analogously). The notation $a \simeq b$ is used if $a \leq cb$ and $b \leq c^\prime a$ with constants $c$ and $c^\prime$ independent of $\kappa$, $\l$, $H_\l$, $L$, and $m$.
\end{remark}
\subsection{$\mathcal{V}$-stable projection}
The $\V$-stable projection introduced in the following, plays a crucial role in the construction of the localized ansatz spaces in Section \ref{sec:pracmethod}. 
The construction of the operator $\Pl$ is based on the quasi-interpolation operator $\Il$, from which it inherits the $\V$-stability. The quasi-interpolation is defined as $\Il:=\mathcal{E}_\l\circ \Pi_\l$, where $\mathcal E_\l\colon \mathbb{Q}_0(\T_\l) \to \mathbb{Q}_1(\T_\l) \cap \V$ ($\mathbb{Q}_1(\T_\l)$ is the space of element-wise tensor-product polynomials of degree up to one) denotes the nodal averaging operator, i.e., for an interior node~$z$, let
\begin{equation*}
(\mathcal E_\l v)(z) 
:= \frac{1}{|T(z)|} \sum\nolimits_{T\in T(z)} v|_T(z)
\end{equation*}
with $T(z):=\{T\in \T_\l: z\in \overline{T}\}$. For boundary nodes, let~$(\mathcal E_\l v)(z):=0$. 
Since $\Pl$ should satisfy $\ker \mathcal{P}_\ell = \ker \Pi_\ell$, we need to correct $\Il v$ such that element averages of $v$ on $\T_\l$ are preserved:
\begin{equation}\label{eq:vstabdef}
\mathcal{P}_\ell v := \mathcal{I}_\ell v +\sum_{T\in \mathcal{T}_\ell} \int_T v -\mathcal{I}_\ell v \dx\, b_T.
\end{equation}
The following stability and approximation properties are shown in  \cite{Peterseim2021} and are based on classical finite element theory \cite{Brenner2008,Ern2017}. For any $v\in \V$
	\begin{align}\label{eq:P}
	\begin{split}
	\|\mathcal{P}_\ell v\|_{L^2(T)} &\lesssim \|v\|_{L^2(N_\l^1(T))},\\  \|\nabla \mathcal{P}_\ell v\|_{L^2(T)} &\lesssim \|\nabla v\|_{L^2(N_\l^1(T))},\\ \|(1-\mathcal{P}_\ell)v\|_{L^2(T)} &\lesssim H_\ell \|\nabla v\|_{L^2(N_\l^1(T))}.
	\end{split}
	\end{align}
\section{Block-diagonal multi-resolution decomposition}\label{sec:idealmethod}

In this section, we introduce a hierarchical  multi-resolution Petrov-Galerkin method in order to discretize \eqref{eq:wfcont}. 
In contrast to \cite{Owhadi2017a,Feischl2020}, where only symmetric and elliptic problems are considered,  we need two different sets of bases for trial and test spaces. This is due to the lack of hermiticity of the Helmholtz problem. 

\begin{figure}[h]
	\begin{tabularx}{\textwidth}{XXXX}
		\includegraphics[width=\linewidth]{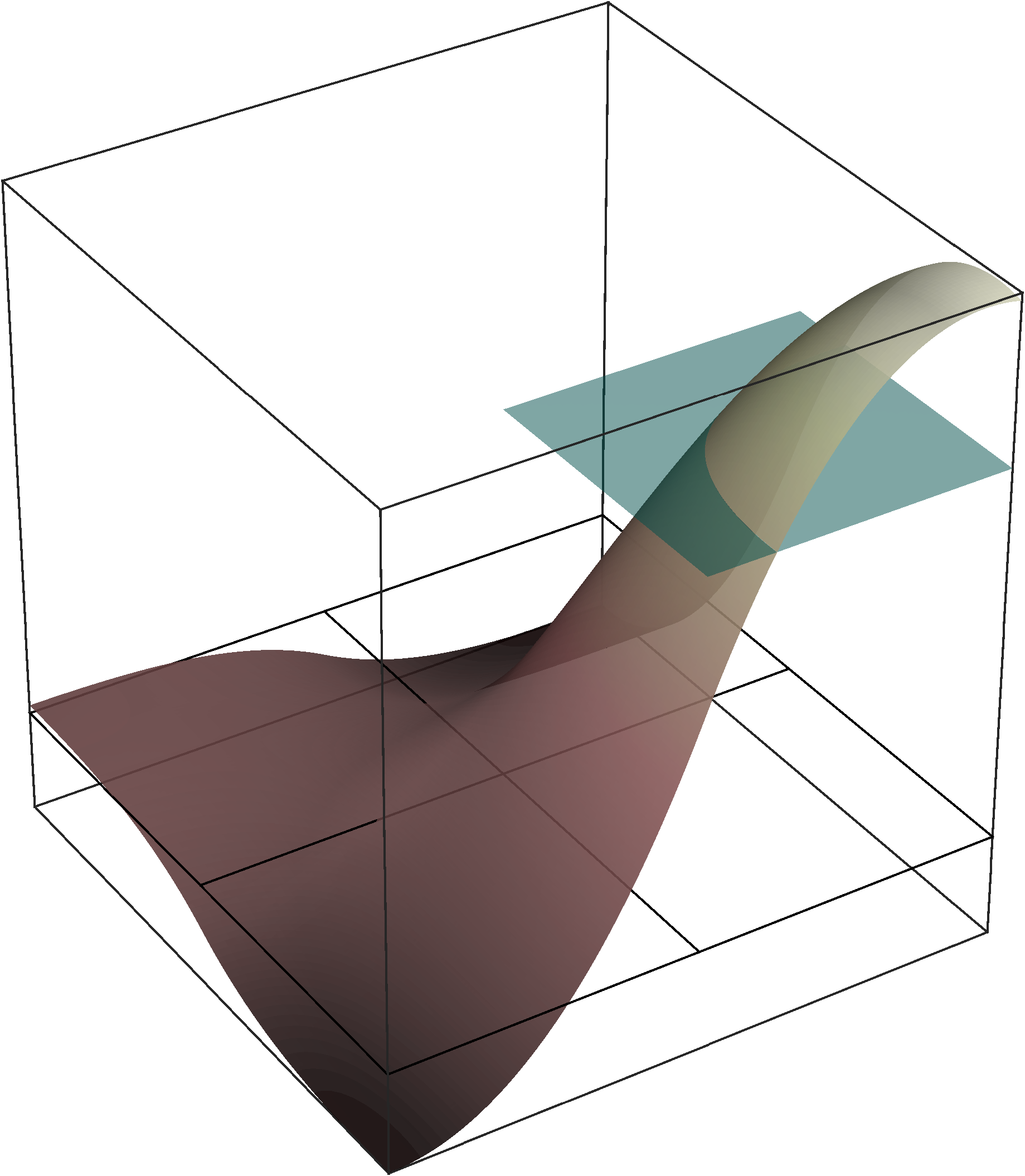}&
		\includegraphics[width=\linewidth]{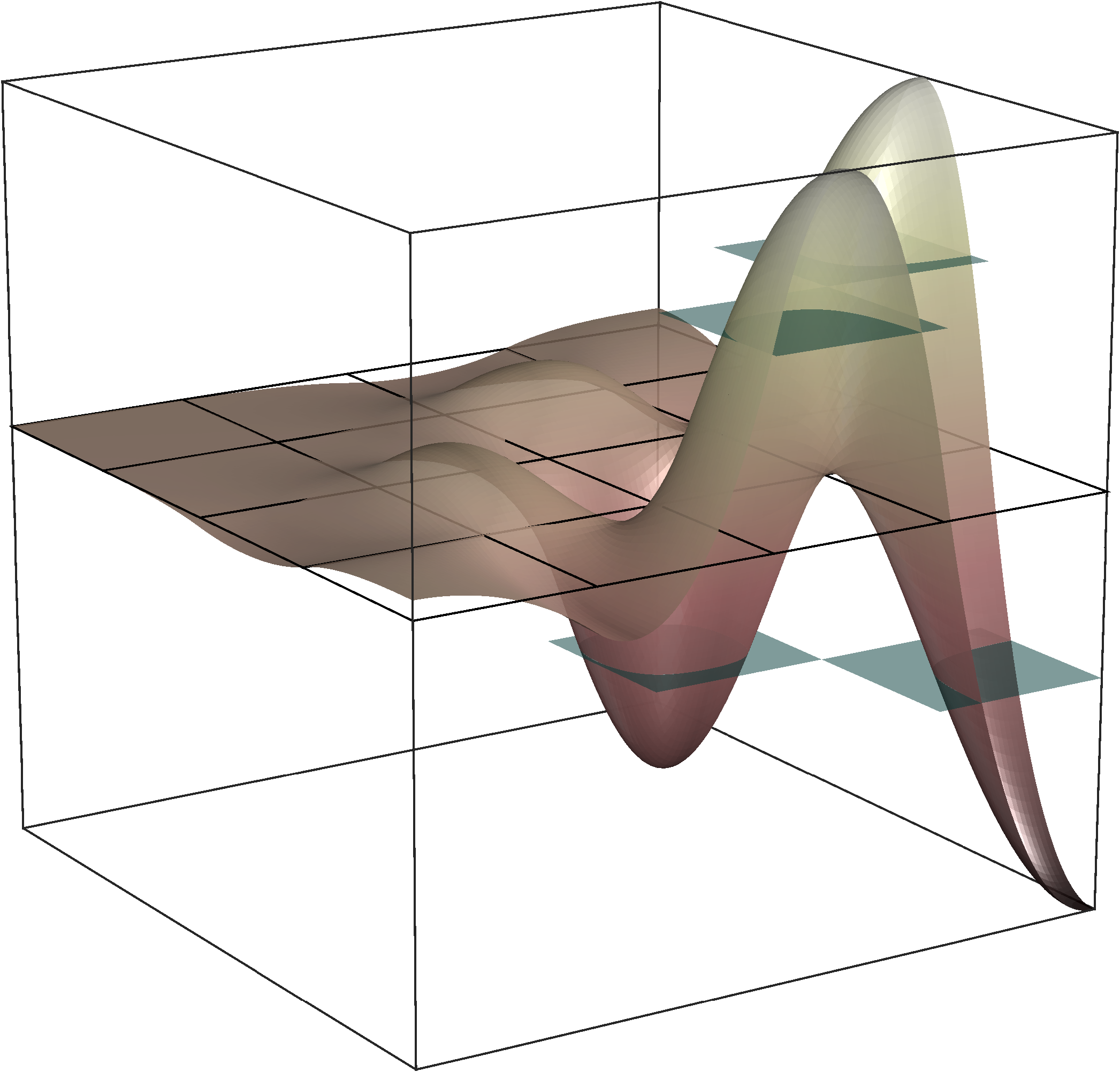}&
		\includegraphics[width=\linewidth]{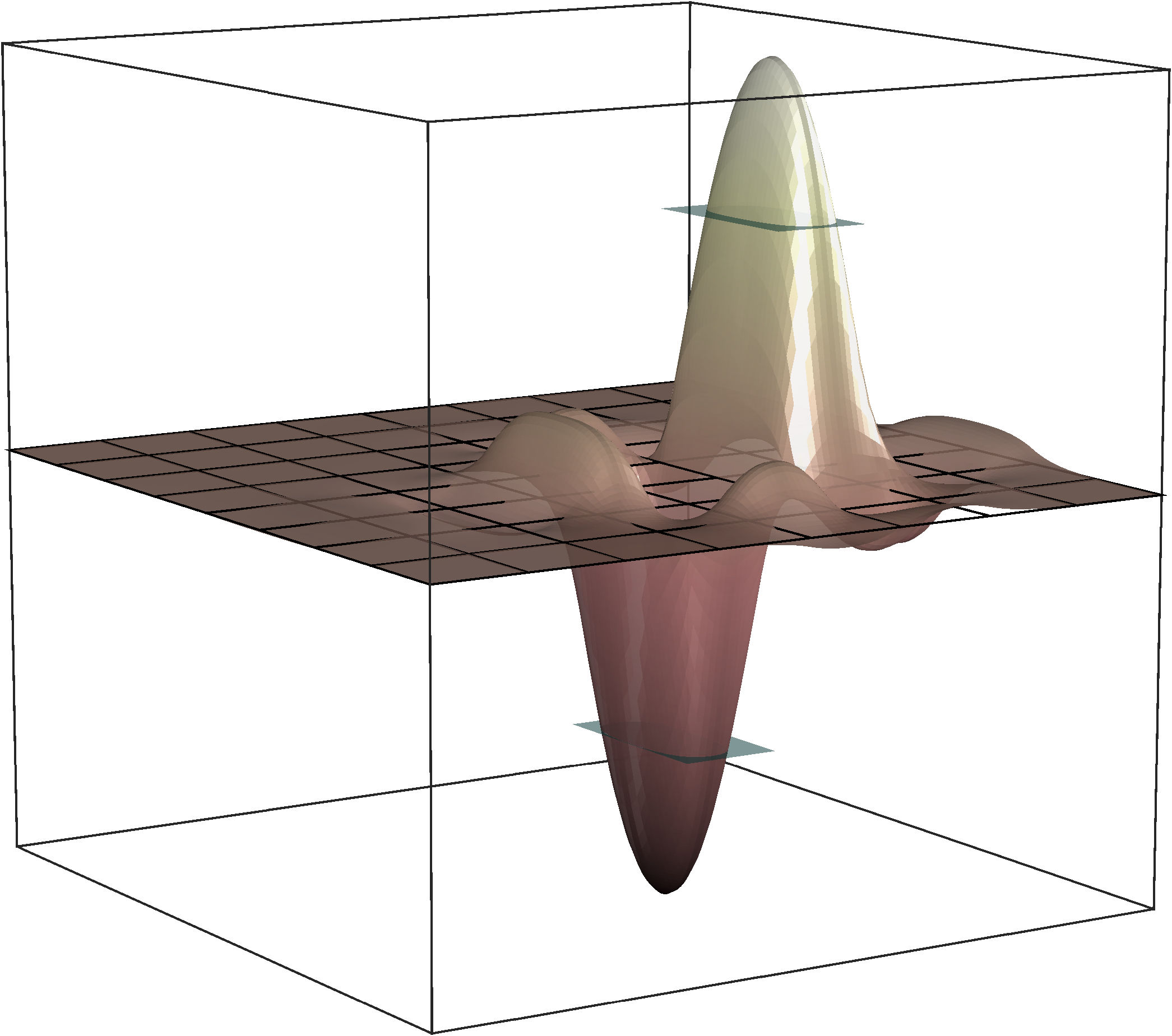}&
		\includegraphics[width=\linewidth]{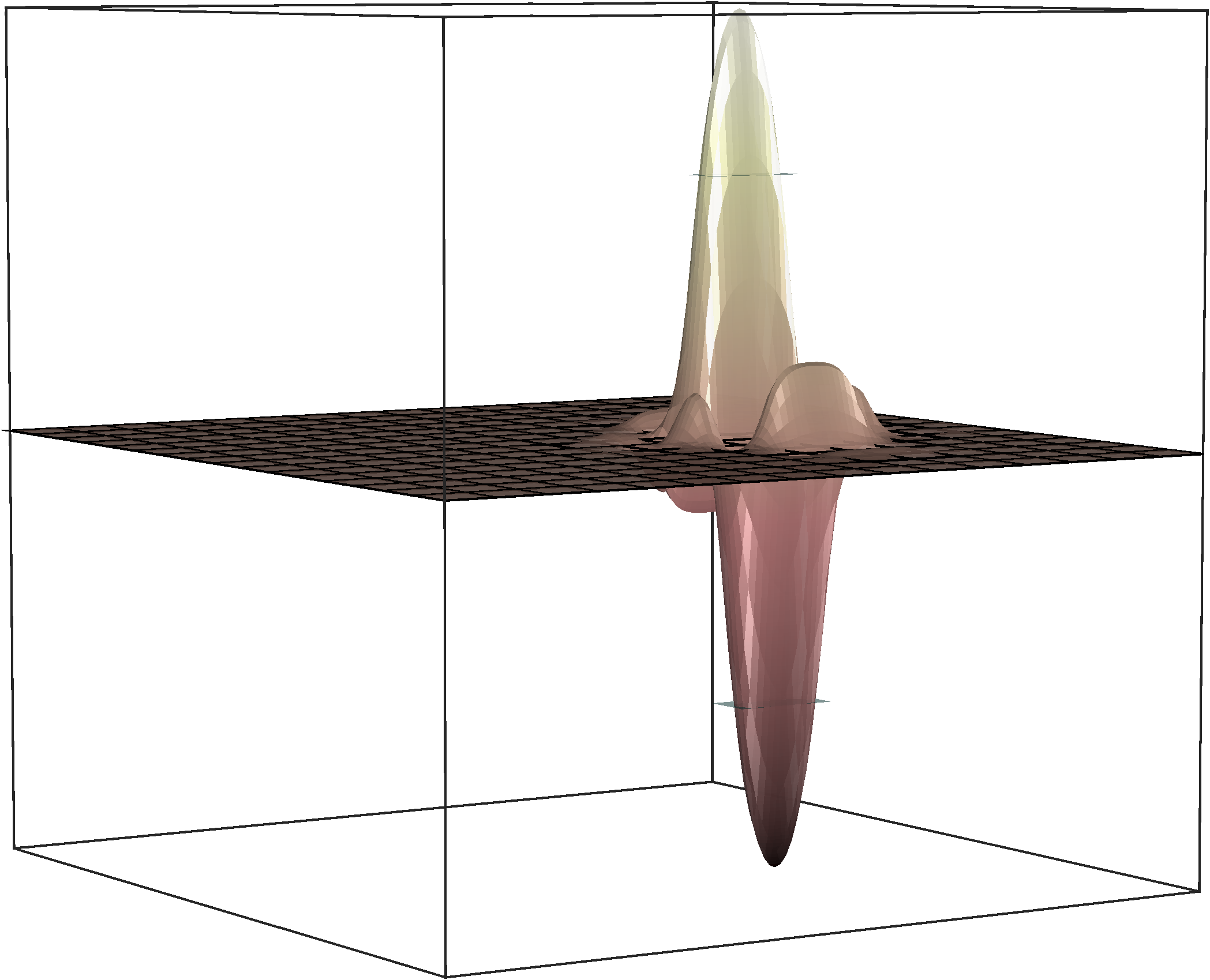}
	\end{tabularx}
	\caption{Let the mesh hierarchy be as for Figure \ref{fig:haarbasis}. The pictures show the real part of one (randomly chosen) basis function  $\bli$ per level. The corresponding Haar basis functions $\phi_{\l,j}$ are indicated in green.}
	\label{fig:idealgambbasis}
\end{figure}

We propose a special choice of bases, which block-diagonalizes the (discretized) operator, see \cite{Feischl2020}. Each basis function is assigned to (exactly) one Haar basis function $\phi_{\l,j}$. Since the basis functions, in general, have global support (c.f. Figure \ref{fig:idealgambbasis}), the ideal method is not suited for practical purposes. A practical method is derived in Section \ref{sec:pracmethod} by localizing the basis functions to element patches. This procedure is justified by the exponential decay of the basis functions (Section \ref{sec:localizetion}).
\subsection{Correction operators}
The concept of correction operators is key in the LOD framework. In the multi-level context, one considers correction operators on every level projecting onto $\W_\l:=  \mathrm{ker}\, \Pi_\l$. The spaces $\W_\l$ can be interpreted as (infinite-dimensional) subspaces of $\V$ consisting of functions oscillating on length scales smaller than $H_\l$. The fine scale spaces are nested, i.e., $\W_1\supset...\supset \W_L$. 

For the correction operators to be well-defined, we need to assume a minimal resolution condition.
{\color{black} 
\begin{assumption}[Resolution condition]\label{assumption:rescond}
	Given the wave number $\kappa$, we assume that the mesh size $H_1$ of $\mathcal T_1$ satisfies 
	$$H_1 \kappa \leq \frac{\pi}{\sqrt{2}}.$$ 
\end{assumption}
}
Due to the lack of hermiticity, we need to define for $\l = 1,...,L$ the two correction operators $\C_\l$ and $\C_\l^*$ as 
\begin{align}\label{eq:defcorr}
\begin{split}
\C_\l: \V\rightarrow \W_\l,\ a(\C_\l v,w) &= a(v,w) \quad \text{for all }w\in \W_\l,\\
\C_\l^*: \V\rightarrow \W_\l,\ a(w,\C_\l^*v) &= a(w,v) \quad \text{for all }w\in \W_\l.
\end{split}
\end{align} 
One can show the identity $\C_\l^*v = \overline{\C_\l\overline{v}}$ for all $v\in \V$. The well-posedness of $\C_\l$ and $\C_\l^*$ follows from the resolution condition.
\begin{lemma}[Well-posedness of corrector problems]\label{lemma:coercivity}
	Under Assumption \ref{assumption:rescond}, $a$ is coercive on $\W_1\times\W_1$, i.e., 
	\begin{equation*}
	\mathfrak{R}a(w,w)\geq \frac{1}{2}  \|\nabla w\|_{L^2(D)}.
	\end{equation*} 
	Additionally, the norms $\|\cdot\|_{\V(D)}$ and $\|\nabla \cdot\|_{L^2(D)}$ are equivalent on $\W_1$ independently of $\kappa$.
	\begin{proof}
		{\color{black}
		$a$ is coercive on $\W_1\times\W_1$, since
		\begin{align*}
		\mathfrak{R}a(w,w) &\geq \tnorm{\nabla w}{D}^2 - \kappa^2\tnorm{w}{D} = \tnorm{\nabla w}{D}^2 - \kappa^2\tnorm{(1-\Pi_1)w}{D}^2\\
		&\geq \tnorm{\nabla w}{D}^2 -  \pi^{-2}(H_1\kappa)^2\tnorm{\nabla w}{D}^2 \geq \frac{1}{2}\tnorm{\nabla w}{D}^2
		\end{align*}
		using Assumption \ref{assumption:rescond}. The norm equivalence of $\tnorm{\nabla \cdot}{D}$ and $\vnorm{\cdot}{D}$ follows from
		\begin{equation*}
		\tnorm{\nabla w}{D}^2\leq \vnorm{w}{D}^2 = \tnorm{\nabla w}{D}^2 + \kappa^2\tnorm{(1-\Pi_1)w}{D}^2\leq \frac{3}{2} \tnorm{\nabla w}{D}^2
		\end{equation*}
		using again the same assumption.}
	\end{proof}
\end{lemma}
\begin{remark}[Well-posedness of correctors]\label{remark:wellposednesscorr}
	The well-posedness of \eqref{eq:defcorr} can be concluded using the Lax--Milgram lemma using the coercivity of $a$ on $\W_1\times\W_1$ shown in Lemma \ref{lemma:coercivity} and the continuity of $a$ \eqref{eq:conta}. In general, it holds that the complementary projection and the projection itself have the same operator norms (c.f. \cite{Szyld2006}). Thus, we get boundedness of $\C_\l$, $\C_\l^*$, $1-\C_\l$, and $1-\C_\l^*$ all with the same constants.
\end{remark}
\subsection{Trial and test spaces}
Using the correction operators, we can define the hierarchical multi-resolution trial and test spaces of the method. For levels $\l = 1,...,L$, let
\begin{align}\label{eq:choicebasis}
\begin{split}
{\widetilde{ \Phi}}_\l &:= \linh\big\{ (1-\C_\l)\{\Pit_\l\phi_{\l,j},j=1,...,N_\l\}\big\},\\ {\widetilde{\Psi}}_\l &:= \linh\big\{ (1-\C^*_\l)\{\Pit_\l\phi_{\l,j},j=1,...,N_\l\}\big\}.
\end{split}
\end{align}
The canonical bases of $\Xtl$ and $\Ytl$ are given by $\{\bli:=(1-\Cl)\Pit_\l\phi_{\l,j},j=1,...,N_\l\}$ and  $\{\blis:=(1-\Cl^*)\Pit_\l\phi_{\l,j},j=1,...,N_\l\}$, respectively.

Let us define the spaces $\Ut := \Xt_1\oplus...\oplus\Xt_L$ and $\Vt := \Yt_1\oplus...\oplus\Yt_L$.
\begin{remark}[Conforming version of $\H_\l$]
	Note that $\Pit_\l\H_\l$ is a conformal version of $\H_\l$. It holds that $\Pit_\l\phi_{\l,j}$ and $\phi_{\l,j}$ have the same element-averages on $\T_\l$. The precise choice of $\Pit_\l$ is not important as long as its image is conforming, it preserves element-averages, and \eqref{eq:Pit} holds.
\end{remark}


\begin{remark}[Relationship to \cite{Peterseim2021}]
	In \cite{Peterseim2021}, the LOD method is introduced using a very general framework. So called 'quantities of interest' determine the precise method. Using the element-averages on the mesh $\T_L$ as quantities of interest, we can prove the relationships
	\begin{equation*}
	(1-\C_L)\V = \Xt_1 \oplus...\oplus\Xt_L\quad\text{ and }\quad(1-\C^*_L)\V = \Yt_1 \oplus...\oplus\Yt_L,
	\end{equation*}
	where  $(1-\C_L)\V$ and $(1-\C_L^*)\V$ are the trial and test spaces from \cite{Peterseim2021}, respectively. This relationship, in general,  is not true for the practical method derived in Section \ref{sec:pracmethod}.
\end{remark}

From Definition \eqref{eq:defcorr}, one can deduce for $k,\l\in\{1,...,L\}$ with $ k\leq \l$ the  'orthogonality' ($a$ is not an inner-product) relations 
\begin{equation}\label{eq:orthogonalitykernel}
a(\Xt_k,\W_\l) = 0 \quad \text{ and }\quad a(\W_\l,\Yt_k) = 0.
\end{equation}

The single-level method discretizing \eqref{eq:wfcont} is based on the trial-test pairing $(\Ut,\Vt)$ and seeks $\tilde{u}\in \Ut$, such that for all $\tilde{v}\in\Vt$
\begin{equation}\label{eq:wfdisc}
a(\tilde{u},\tilde{v}) = \tsp{f}{\tilde{v}}{D}.
\end{equation}
Henceforth, we refer to \eqref{eq:wfdisc} as ideal method, since this method only has theoretical purposes and cannot be applied in practice; see Section \ref{sec:pracmethod}.

\subsection{Stability and accuracy of ideal method}
Since  \eqref{eq:wfdisc} coincides with the ideal method from \cite{Peterseim2017a}, also the stability result and convergence results are the same.

\begin{lemma}[Stability of ideal method]\label{lemma:infsupUtVt}
	Let Assumption \ref{assumption:rescond} be satisfied. Then the trial space $\Ut$ and the test space $\Vt$ satisfy the inf-sup condition
	\begin{equation*}
	0<\alpha^{(2)}\leq \inf_{\tilde u\in \Ut}\sup_{\tilde v\in\Vt} \frac{\mathfrak{R}a(\tilde{u},\tilde{v})}{\vnorm{\tilde{u}}{D}\vnorm{\tilde{v}}{D}}
	\end{equation*}
	with $\alpha^{(2)}\simeq ( c_{\mathrm{stab}}(\kappa)\kappa)^{-1}$.
	\begin{proof}
		For the proof, see Appendix \ref*{sec:appendix}.
	\end{proof}
\end{lemma}

\begin{lemma}[Accuracy of ideal method]\label{lemma:errideal}
	{\color{black}
	Let $u\in\V$ be the solution of \eqref{eq:wfcont} for the right-hand side $f$. If  Assumption \ref{assumption:rescond} is satisfied, then $\tilde u = (1-\C_L )u\in \Ut$ is the unique solution of \eqref{eq:wfdisc}, that is, the Petrov-Galerkin approximation of $u$ in the subspace $\Ut$ with respect to the test space $\Vt$.  Moreover, there is a constant $c$ independent of $\kappa, H_L$ and $L$ such that for any $f\in H^s(D)$, $s\in[0,1]$
	\begin{equation*}
	\vnorm{u-\tilde{u}}{D} \leq c H_L^{1+s} \|f\|_{H^s(D)}.
	\end{equation*}
	}
	\begin{proof}
		For the proof, see Appendix \ref{sec:appendix}.
	\end{proof}
\end{lemma}
\subsection{Block-diagonalization}
The special choice of the hierarchy of bases \eqref{eq:choicebasis} decouples the equations on every level; see also \cite{Owhadi2017a,Feischl2020}.
\begin{lemma}['Orthogonality' of levels]\label{lemma:orth}
	Let Assumption \ref{assumption:rescond} be satisfied. Consider $k,\l \in \{1,...,L\}$, with $k\neq \l$. Then it holds for all $\xtk\in \Xt_k,\ytl\in \Yt_\l$
	\begin{equation*}
	a(\xtk,\ytl) = 0.
	\end{equation*}
	\begin{proof}
		First suppose that $k>\l$. We show that $\Xt_k \subset \W_\l$. Let $\Xt_k \ni\xtk = (1-\C_k)\Pit_k \phi_k$ with $\phi_k \in \linh \{\H_k\}$ be arbitrary, then 
		\begin{equation*}
		\Pi_\l(1-\C_k)\Pit_k \phi_k = \Pi_\l\Pi_k(1-\C_k)\Pit_k \phi_k = \Pi_\l\phi_k= 0.
		\end{equation*}
		Together with \eqref{eq:orthogonalitykernel}, this implies the assertion. If $k<\l$, one can similarly show that $\Ytl \subset \W_k$, which again yields the assertion. 
	\end{proof}
\end{lemma}
\begin{corollary}[]
	Let Assumption \ref{assumption:rescond} be satisfied, then for $\l\in\{2,...,L\}$ it holds that
	\begin{equation*}
	\Xt_\l\subset \W_{\l-1}\subset \W_1,\quad \text{ and }\quad \Yt_\l\subset \W_{\l-1}\subset \W_1.
	\end{equation*}
	\begin{proof}
		This is an immediate consequence of the proof of Lemma \ref{lemma:orth}.
	\end{proof}
\end{corollary}
Using Lemma \ref{lemma:orth}, we can rewrite \eqref{eq:wfdisc} into decoupled Petrov-Galerkin problems on every level, which 
seek $\xtl\in\Xt_\l$, such that for all $\ytl\in\Yt_\l$ 
\begin{equation}\label{eq:wfdiscml}
a(\xtl,\ytl) = (f,\ytl)_{L^2(D)}.
\end{equation}
The solution $\tilde{u}\in\Ut$ of \eqref{eq:wfdisc} is then obtained by
\begin{equation}\label{eq:wfdiscmladdup}
\tilde{u} = \tilde{\varphi}_1+...+\tilde{\varphi}_L.
\end{equation}

\begin{remark}[Variants of the method]
	Due to \eqref{eq:orthogonalitykernel}, we obtain for $\xtl = (1-\Cl)\Pit_\l\varphi_\l$,  $\ytl = (1-\Cl^*)\Pit_\l\psi_\l$ with $\varphi_\l,\psi_\l\in\linh\{\H_\l\}$ 
	$$a(\xtl,\ytl) =
	 a(\Pit_\l\varphi_\l,(1-\Cl^*)\Pit_\l\psi_\l) = a((1-\Cl)\Pit_\l\varphi_\l,\Pit_\l\psi_\l).$$
	This motivates Petrov-Galerkin methods for \eqref{eq:wfdiscml} which are based on the trial-test pairings $(\linh\{\Pit_\l\H_\l\},\Ytl)$ and $(\Xtl,\linh\{\Pit_\l\H_\l\})$; see also \cite{Gallistl2015,Peterseim2016a}.
\end{remark}
The following lemma shows that the sub-scale problems \eqref{eq:wfdiscml} are well-posed. 
\begin{lemma}[Stability of sub-scale problems for $\l\geq 2$]\label{lemma:infsupXtlYtl}
	Let Assumption \ref{assumption:rescond} be satisfied. Then, for $\l\geq 2$,  the trial space $\Xt_\l$ and the test space $\Yt_\l$ satisfy the inf-sup condition
	\begin{equation*}
	0<\alpha^{(3)}_\l\leq\inf_{\xtl\in \Xt_\l}\sup_{\ytl\in\Yt_\l} \frac{\mathfrak{R}a(\xtl,\ytl)}{\vnorm{\xtl}{D}\vnormf{\ytl}{D}}
	\end{equation*}
	with $\alpha^{(3)}_\l \simeq 1$.
	\begin{proof}
		All $\xtl\in\Xtl$ can be written as $\xtl = (1-\C_\l)\Pit_\l\phi_\l$ with $\phi_\l\in\linh\{\H_\l\}$. Defining $\ytl:=(1-\C_\l^*)\Pit_\l\phi_\l$, we obtain
		\begin{align*}
		a(\xtl,\ytl) = a((1-\C_\l)\Pit_\l \phi_\l ,(\C_\l^*-\C_\l)\Pit_\l\phi_\l) + a(\xtl,\xtl) = a(\xtl,\xtl).
		\end{align*} 
		Using that $\xtl\in\W_1$ for $\l\geq 2$, yields
		\begin{align*}
		\mathfrak{R} a(\xtl,\ytl) = \mathfrak{R}a(\xtl,\xtl) \gtrsim \tnorm{\nabla \xtl}{D}^2\gtrsim  \vnorm{\xtl}{D}^2,
		\end{align*}
		where we used Assumption \ref{assumption:rescond} and Lemma \ref{lemma:coercivity}. Since $\ytl\in \W_1$, we obtain 
		\begin{align*}
		\vnormf{\ytl}{D}^2&\lesssim \mathfrak{R}a(\ytl,\ytl) =\mathfrak{R}a((\C_\l-\C_\l^*)\Pit_\l\phi_\l,(1-\C_\l^*)\Pit_\l\phi_\l)+\mathfrak{R}a(\xtl,\ytl)\\
		&= \mathfrak{R}a(\xtl,\ytl)\lesssim \vnorm{\xtl}{D}\vnormf{\ytl}{D}.
		\end{align*}
		Combining the previous estimates yields the assertion.
	\end{proof}
\end{lemma}
\begin{remark}[Stability of sub-scale problem for $\l =1$]\label{infsupl1}
	Since $\Xtl\nsubseteq\W_1$ for $\l = 1$, the argument from Lemma \ref{lemma:infsupXtlYtl} is not valid in this case. For $\l = 1$, one obtains a stability estimate similarly to Lemma \ref{lemma:infsupUtVt} with inf-sup constant $\alpha^{(3)}_1 \simeq ( c_{\mathrm{stab}}(\kappa)\kappa)^{-1}$.
\end{remark}

\section{Localized numerical correctors}\label{sec:localizetion}
In this section, we decompose the correction operator into element-correctors, which decay exponentially fast away from the corresponding element (Lemma \ref{lemma:expdec}). The decay makes it possible to localize the element-corrector problems to element patches. Approximation properties are shown (Lemma \ref{lemma:locerrcorr}). 
\subsection{Exponential decay of element-correctors}
For $\l\in\{1,...,L\}$ and all $T\in\T_\l$,  we define element-correctors $\Clt, \Clt^*$ for all $v\in\V$ as 
\begin{align}\label{eq:elemcorrdef}
\Clt: \V\rightarrow \W_\l,\ a(\Clt v,w) &= a_T(v,w) \quad \text{for all }w\in \W_\l
\end{align} 
and $\Clt^* v:=\overline{\Clt \overline{v}}$. Here, we used the restricted sesquilinear form $a_T$ defined as 
\begin{equation*}
a_T(v,w) := \tsp{\nabla v}{\nabla w}{T} - \kappa^2\tsp{v}{w}{T} -i\kappa \tsp{v}{w}{\partial T\cap \partial D}.
\end{equation*}
\begin{remark}[Well-posedness of element-correctors]
	The element-correctors are well-posed due to the same arguments as in Remark \ref{remark:wellposednesscorr}. Similarly as before, we get boundedness of $\Clt$, $\Clt^*$, $1-\Clt$, and $1-\Clt^*$ all with the same constants.
%
\end{remark}


Next, we show that the moduli of the element-correctors decay exponentially fast.

\begin{lemma}[Exponential decay of element-correctors]\label{lemma:expdec}
	{\color{black}
	If Assumption \ref{assumption:rescond} is satisfied, then there exists a constant $0<\beta<1$ independent of $\kappa$, $H_\l$ and $L$ such that for all $\l\in\{1,...,L\}$, $T\in\T_\l$, $v\in \V$ and $m \in\mathbb{N}$, the element correctors $\Clt v$ satisfy}
	\begin{equation*}
	\tnorm{\nabla \Clt v}{D\backslash N_\l^m(T)}\leq \beta^m\tnorm{\nabla \Clt v}{D}.
	\end{equation*}
	An analogous decay result holds for $\Clt^*$.
	\begin{proof}
		For the proof, see Appendix \ref{sec:appendix}.
	\end{proof}
\end{lemma}
{\color{black}
\begin{remark}[Rate of decay]
	From the proof of Lemma \ref{lemma:expdec}, it becomes clear that, if Assumption \ref{assumption:rescond} is satisfied, the rate of decay $\beta$ is independent of $\kappa$, $H_\ell$ and $L$.
\end{remark}}

\subsection{Localized correctors}
The exponential decay motivates the localization of the \linebreak (global) problems \eqref{eq:elemcorrdef} to patches, i.e., for $\l\in\{1,...,L\}$ and an oversampling parameter $m\in \mathbb{N}$, we substitute the (global) ansatz space $\W_\l$ by its localized counterparts
\begin{equation*}
\W_\l^m(T):=\{w\in \W_\l\,:\, \supp\,w\subset N^m_\l(T)\}\subset \W_\l.
\end{equation*}
We then define the localized element-correctors for all $v\in\V$ as the solution of
\begin{align}\label{eq:locelemcorr}
\Cltm: \V\rightarrow \W_\l^m(T),\ a(\Cltm v,w) &= a_T(v,w) \quad \text{for all }w\in \W_\l^m
\end{align} 
and set $\Cltms v:=\overline{\Cltm \overline{v}}$. The element-correctors $\Clt$ and $\Clt^*$ are constructed such that their sum (over all $T\in\T_\l$) equals the global correctors $\Cl$ and $\Cl^*$, respectively. Hence, the sum of the localized element-correctors
\begin{equation*}
\Clm := \sum_{T\in \mathcal{T}_\ell}\Cltm\quad \text{ and }\quad \Clms := \sum_{T\in \mathcal{T}_\ell}\Cltms
\end{equation*}
should be good approximation to the actual corrector $\Cl$ and $\Cl^*$, respectively. The (exponential) approximation properties are shown in the next lemma.
\begin{lemma}[Localization error of correctors]\label{lemma:locerrcorr}
	If Assumption \ref{assumption:rescond} is satisfied, then, 
	{\color{black}there exists a constant $c_m$ depending polynomially on $m$, but independent of $\kappa$, $H_\ell$ and $L$ such that for all $\l \in \{1,\dots,L\}$, $v\in\V$ and $m\in\mathbb{N}$, it holds
	\begin{equation*}
	\tnorm{\nabla(\Clm-\Cl)v}{D}\leq c_m \beta^m \vnorm{v}{D},
	\end{equation*}
	where $\beta$ is the constant from Lemma \ref{lemma:expdec}.} An analogous result holds for $\Clms-\Cl^*$.
	\begin{proof}
		For the proof, see Appendix \ref{sec:appendix}.
	\end{proof}
\end{lemma}

\section{Sparsification}\label{sec:pracmethod}
\begin{figure}[h]
		\includegraphics[height=.4\linewidth]{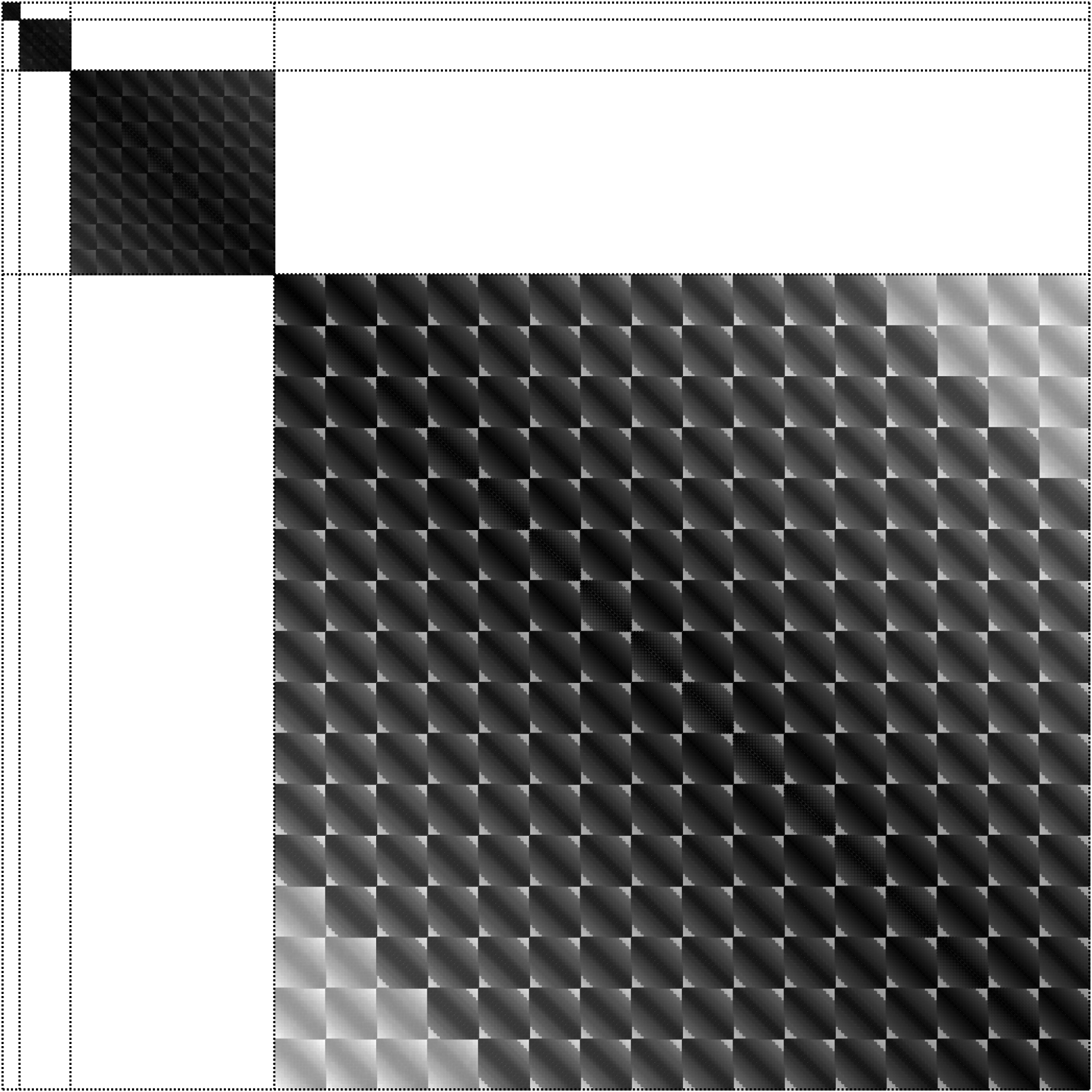}
		\hspace{1cm}
		\includegraphics[height=.4\linewidth]{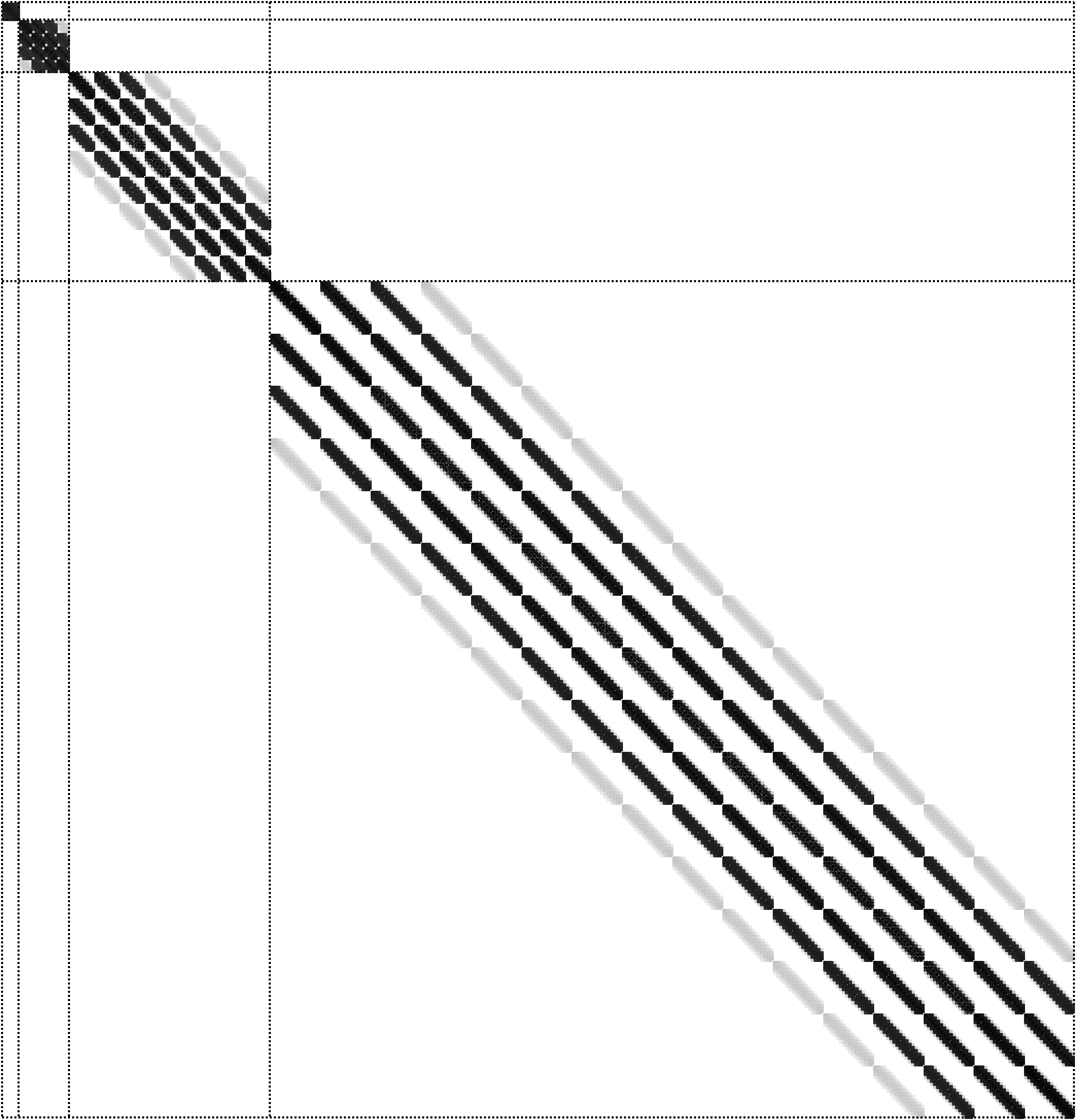}
\caption{Schematic demonstration of the sparsification of the block-diagonal system matrix. Left: Block-matrix corresponding to  ideal method \eqref{eq:wfdiscml}; right: Block-matrix corresponding to proposed practical multi-resolution method \eqref{eq:wfdiscmlloc}. We use a logarithmic gray-scale as color-coding, where dark gray indicates a large modulo of the corresponding entry; zero entries are colored white.}
\label{fig:sparsification}
\end{figure}

In this section, we introduce localized variants of the hierarchical multi-resolution trial and test spaces from Section \ref{sec:idealmethod}. The localization procedure results in a practical method with sparse block-diagonal system matrix; see Figure \ref{fig:sparsification} for a comparison of the sparsity patterns for the ideal and localized methods.

\subsection{Localized trial and test spaces}
Here, we introduce the novel localization strategy that improves stability even for the elliptic case. Instead of just replacing the correctors in \eqref{eq:choicebasis} by their localized counterparts, we additionally utilize the operator $\Pl$ in the definition of the localized trial and test spaces. For $\l = 1,...,L$, let 
\begin{align}\label{eq:deflocttsp}
\begin{split}
\Xtlm &:= \linh\big\{(1-\Clm)\{\Pl\Pit_\l\phi_{\l,j},j=1,...,N_\l\}\big\},
\\ \Ytlm &:= \linh\big\{(1-\Clms)\{\Pl\Pit_\l\phi_{\l,j},j=1,...,N_\l\}\big\}.
\end{split}
\end{align}
The canonical bases of $\Xtlm$ and $\Ytlm$ are given by
\begin{align}\label{eq:basesXtlmYtlm}
\begin{split}
\big\{\blim&:=(1-\Clm)\Pl \Pit_\l \phi_{\ell,j},j=1,...,N_\l\big\},\\
\big\{\blims&:=(1-\Clms)\Pl \Pit_\l \phi_{\ell,j},j=1,...,N_\l\big\}.
\end{split}
\end{align} 
\begin{figure}[h]
	\begin{tabularx}{\textwidth}{XXXX}
		\includegraphics[width=\linewidth]{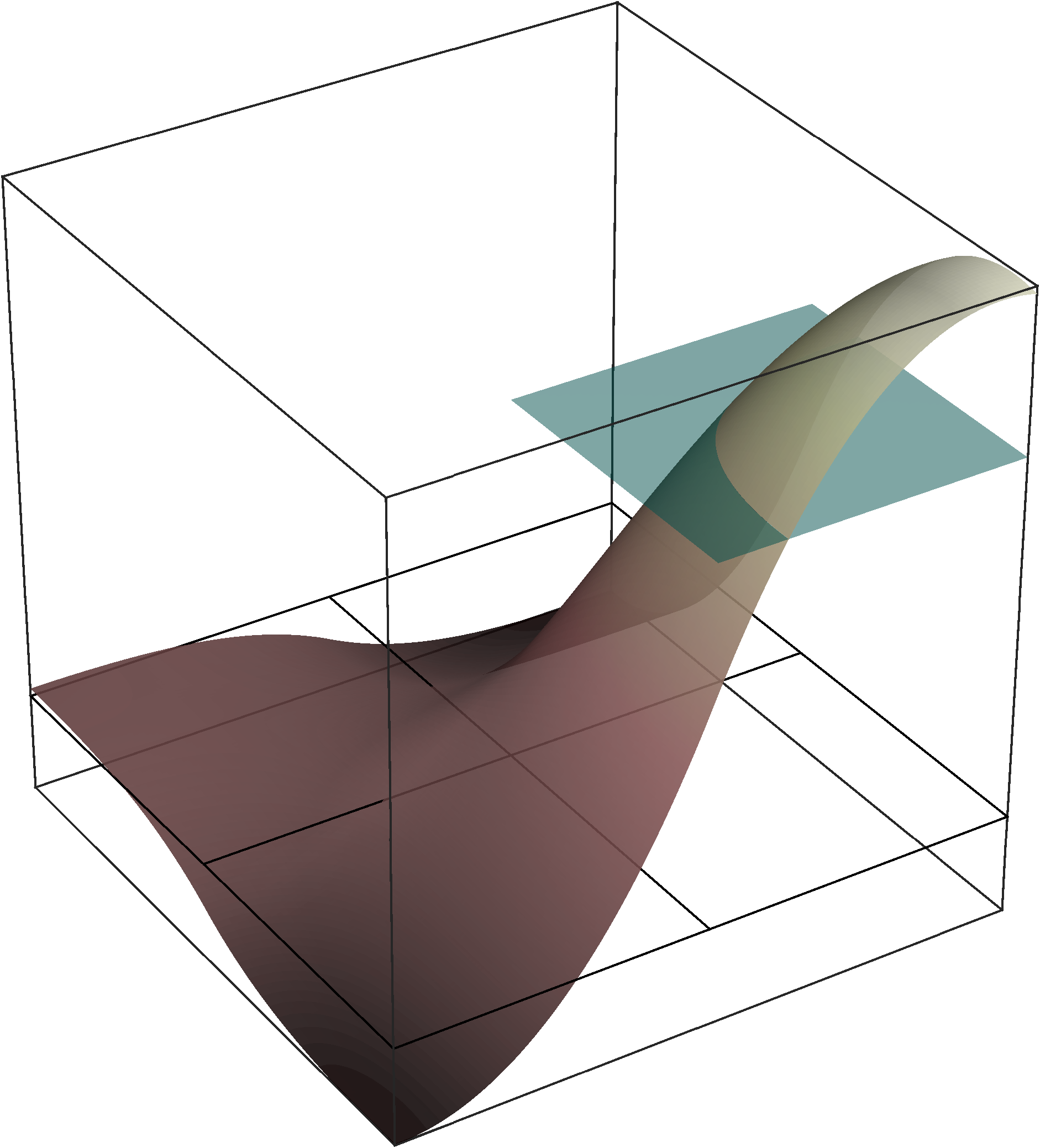}&
		\includegraphics[width=\linewidth]{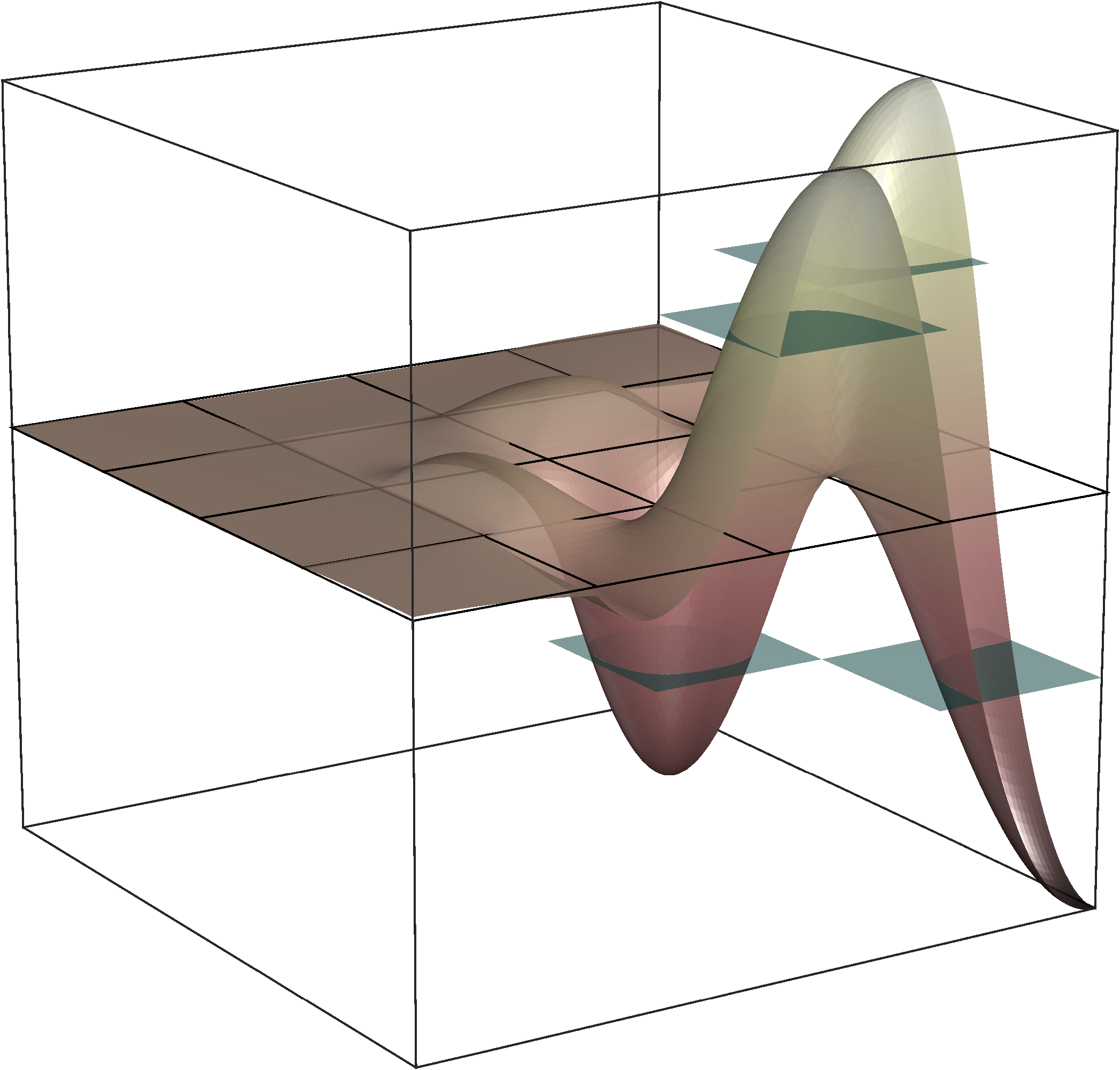}&
		\includegraphics[width=\linewidth]{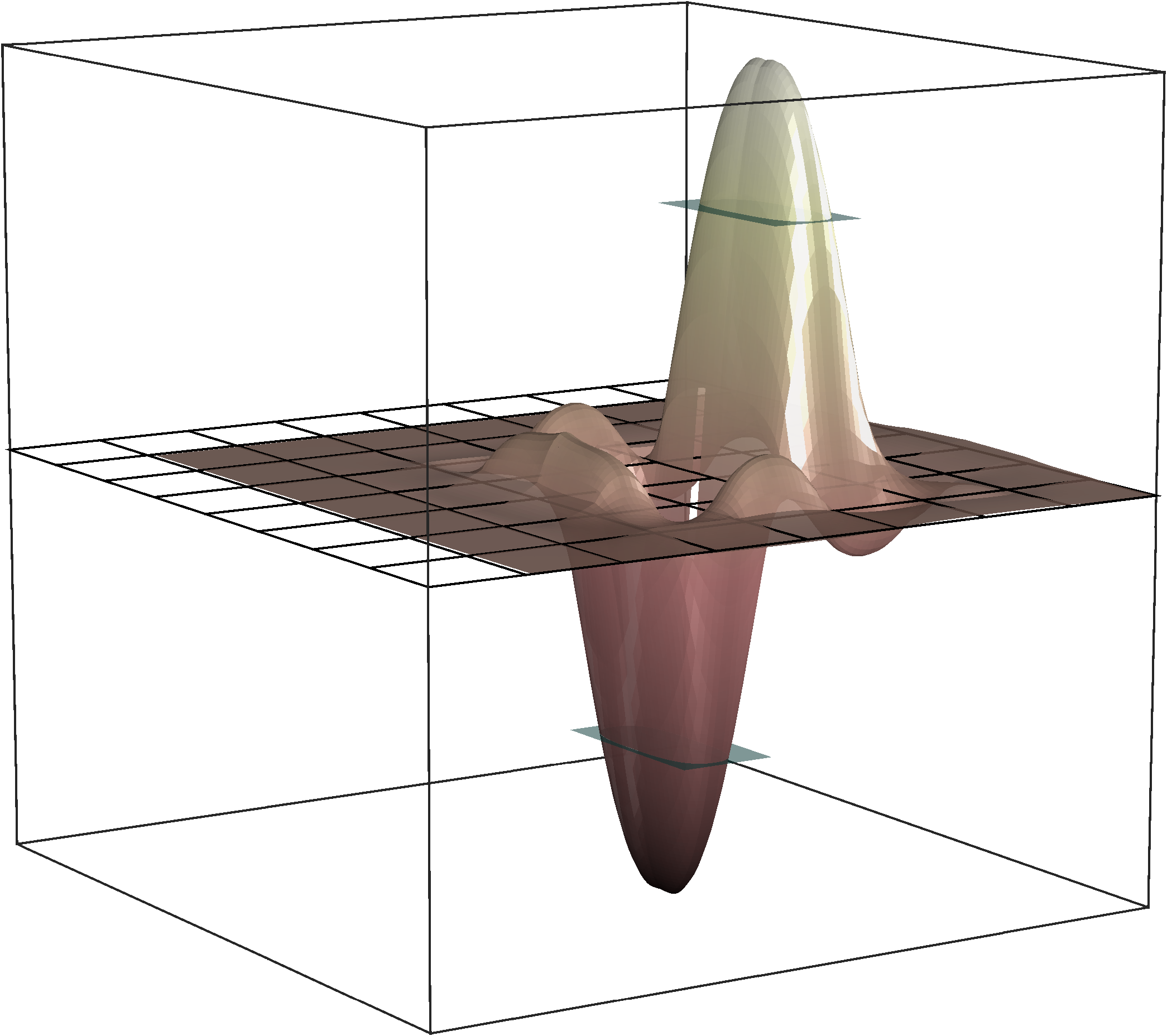}&
		\includegraphics[width=\linewidth]{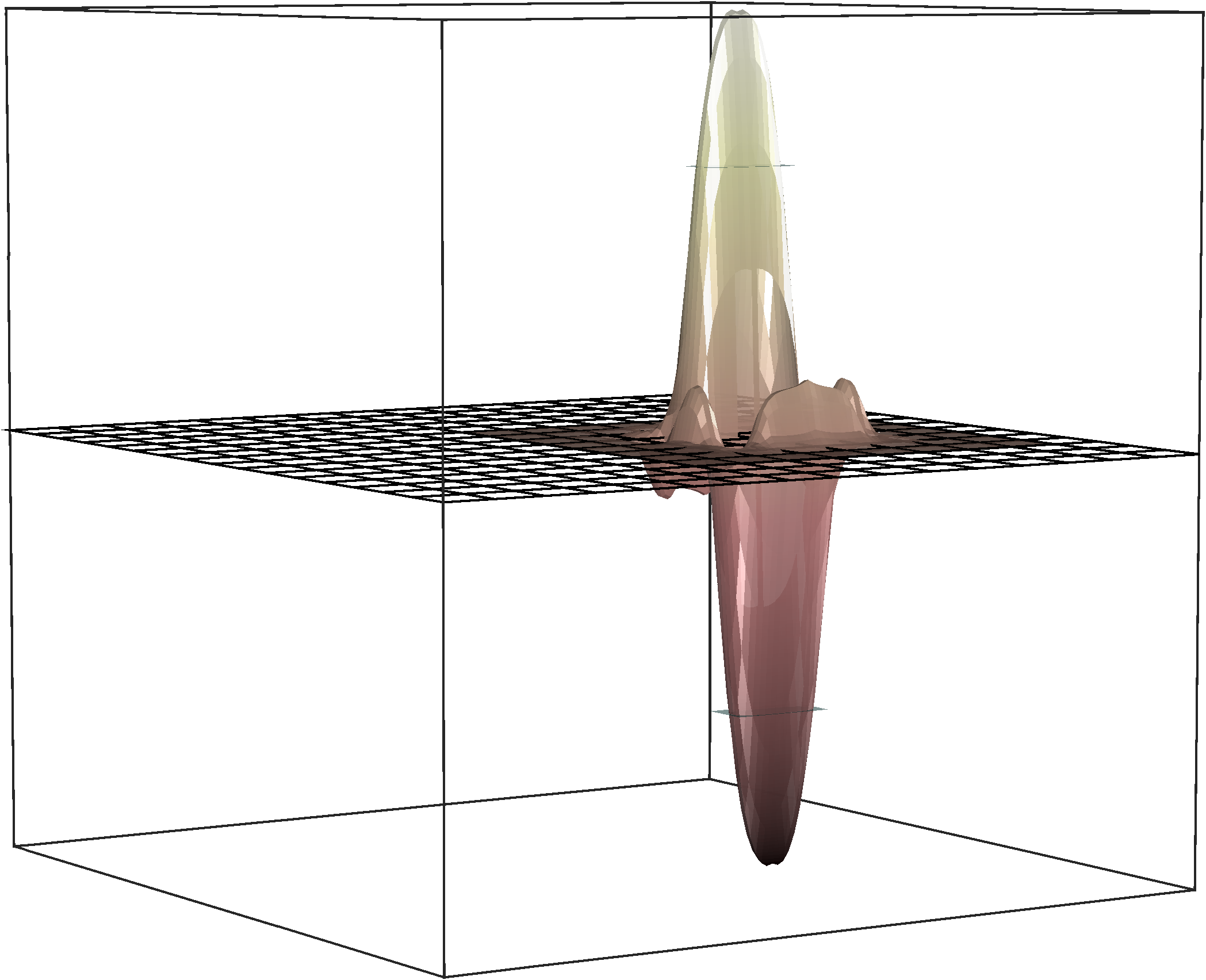}
	\end{tabularx}
	\caption{Let the mesh hierarchy be as for Figure \ref{fig:haarbasis}. The pictures show the real part of one (randomly chosen) basis function  $\bli^1$ per level. The corresponding Haar basis functions $\phi_{\l,j}$ are indicated in green.}
	\label{fig:gambbasis}
\end{figure}

\begin{remark}[Effect of operator $\mathcal P_\ell$]\label{remark:stabH1}
	Since $\Pl$ preserves element-averages, we have $(1-\C_l)\Pit\phi_\l = (1-\C_\l)\Pl\Pit_\l\phi_\l$ for all $\phi_\l\in\linh\{\H_\l\}$ (analogously for $\Cl^*$). Thus, \eqref{eq:deflocttsp} is consistent with \eqref{eq:choicebasis} in the ideal (non-localized) case. The trial and test spaces \eqref{eq:deflocttsp} define  a method that is uniformly stable in $H_\l$. In contrast, omitting  $\Pl$ in \eqref{eq:deflocttsp}  would result in a method, where the oversampling parameter $m$ had to be increased if $H_1$ was decreased; see \cite{Owhadi2017a,Owhadi2019,Maier2020a}.
\end{remark}
{\color{black}
\begin{remark}[Optimality of exponential decay of basis functions]
	In the recent work \cite{Hauck2021}, a multi-scale basis is constructed with super-exponential localization properties. While some theoretical aspects are still unclear, numerical experiments show the supremacy of the novel localization approach compared to existing approaches. 
\end{remark}
}
\begin{remark}[Exponential smallness of off-diagonal blocks]\label{remark:notbd}
	After localization, the operator is not exactly block-diagonal anymore. However, it can be shown that the off-diagonal blocks are exponentially small; see proof of Theorem \ref{lemma:infsuploc}. Thus, for sufficiently large $m$, this error does not impact stability and convergence results. 
\end{remark}

Let us propose the practical variant of the multi-level method \eqref{eq:wfdiscml}, \eqref{eq:wfdiscmladdup}: First, on every level, solve Petrov-Galerkin problems, which seek $\xtlm\in\Xtlm$, such that for all $\ytlm\in\Ytlm$ 
\begin{equation}\label{eq:wfdiscmlloc}
a(\xtlm,\ytlm) = (f,\ytlm)_{L^2(D)}.
\end{equation}
Second, define the solution as the sum of the level contributions
\begin{equation}\label{eq:wfdiscmladduploc}
\tilde{u}^m := \tilde{\varphi}_1^m+...+\tilde{\varphi}_L^m.
\end{equation}
\subsection{Stability and accuracy of the practical multi-level method}
Provided that the localized method is sufficiently close to the ideal method, the stability and convergence results of the ideal method (Lemma \ref{lemma:infsupUtVt} and \ref{lemma:errideal}) carry over to the localized method. The closeness is ensured by a coupling of the oversampling parameter $m$ to the stability constant $c_{\mathrm{stab}}(\kappa)$. We assume a polynomial dependence of the stability constant on the wave number $\kappa$; for a discussion of this assumption, see Section \ref{sec:Helmhprop}.

\begin{assumption}[Polynomial-in-$\kappa$ stability]\label{assumption:poly}
	We assume that $c_{\mathrm{stab}}$ depends polynomially on $\kappa$, i.e., 
	\begin{equation*}
	c_{\mathrm{stab}}(\kappa) \lesssim \kappa^n.
	\end{equation*}
\end{assumption}
\begin{theorem}[Stability of localized sub-scale problems]\label{lemma:infsuploc}
		Let Assumption \ref{assumption:rescond} and \ref{assumption:poly} be satisfied, as well as the oversampling condition
		\begin{equation}\label{eq:oversampling1}
		m\gtrsim|\log(\beta)|^{-1}
		|\log(\kappa^{n+1}c_m)|
		\end{equation}
		with $\beta$ and $c_m$ from Lemma \ref{lemma:locerrcorr}. Then, for  $\l 1,...,L$, the trial space $\Xtlm$ and the test space $\Ytlm$ satisfy the inf-sup condition
	\begin{equation*}
	0<\alpha^{(4)}_\l\leq\inf_{\xtlm\in \Xtlm}\sup_{\ytlm\in\Ytlm} \frac{\mathfrak{R}a(\xtlm,\ytlm)}{\vnormf{\xtlm}{D}\vnormf{\ytlm}{D}},
	\end{equation*}
	where $\alpha^{(4)}_\l \simeq \alpha^{(3)}_\l$, with constants $\alpha^{(3)}_\l$ from Remark \ref{infsupl1} $(\l = 1)$ and Lemma \ref{lemma:infsupXtlYtl} $(\l\geq 2)$.
	\begin{proof}
	Let $\l\in\{1,...,L\}$ be fixed and let $\xtlm\in\Xtlm$ be arbitrary. Define $\xtl:=(1-\C_\l)\Pl\xtlm$ and $\ytlm:=(1-\Clms)\Pl\ytl$, where $\ytl\in\Ytl$ is chosen such that (Lemma \ref{lemma:infsupXtlYtl})
	\begin{equation*}
	\mathfrak{R}a(\xtl,\ytl)\geq \alpha^{(3)}_\l\vnorm{\xtl}{D}\vnormf{\ytl}{D}.
	\end{equation*}
	We obtain
	\begin{align*}
	\mathfrak{R}a(\xtlm,\ytlm) = \mathfrak{R} a(\xtlm,\ytl)+\mathfrak{R}a(\xtlm,\ytlm-\ytl) = \mathfrak{R} a(\xtl,\ytl)+\mathfrak{R}a(\xtlm,\ytlm-\ytl).
	\end{align*}
	The second term can be bounded by
	\begin{align*}
	\big\vert\mathfrak{R}a(\xtlm,\ytlm-\ytl)\big\vert&\lesssim \vnormf{\xtlm}{D}\vnormf{\ytlm-\ytl}{D} = \vnormf{\xtlm}{D}\vnormf{(\Clms-\C_\l^*)\Pl\ytl}{D}\\
	&\lesssim c_m\beta^m\vnormf{\xtlm}{D}\vnormf{\ytl}{D}.
	\end{align*}
	It remains to estimate the $\V$-norm of $\xtl$ and $\ytl$ against the $\V$-norm of its localized counterparts and vice-versa. With \eqref{eq:P} and Remark \ref{remark:wellposednesscorr}, we obtain
	\begin{align*}
	\vnormf{\xtl}{D} = \vnormf{(1-\C_\l)\Pl\xtlm}{D}\lesssim \vnormf{\Pl\xtlm}{D}\lesssim \vnormf{\xtlm}{D},
	\end{align*}
	Furthermore, with Lemma \ref{lemma:locerrcorr}, we get 
	\begin{align*}
	\vnormf{\xtlm}{D} &= \vnormf{(1-\Clm)\Pl\xtlm}{D}\leq \vnormf{(1-\Cl)\Pl \xtlm}{D}+\vnormf{(\Clm-\C_\l)\Pl\xtlm}{D}\\
	&\lesssim \vnormf{\xtl}{D} + c_m \beta^m\vnormf{\Pl\xtlm}{D}\lesssim (1+c_m\beta^m)\vnormf{\xtl}{D}\lesssim \vnormf{\xtl}{D},
	\end{align*}
	where we used that $c_m\beta^m$ can be bounded independently of $m$.
	The corresponding estimates for $\ytl$ and $\ytlm$ can be obtained using similar arguments. 

	We conclude 
	\begin{align*}
	\mathfrak{R}a(\xtlm,\ytlm) &= \mathfrak{R} a(\xtl,\ytl)+\mathfrak{R}a(\xtlm,\ytlm-\ytl)\\
	&\gtrsim \alpha^{(3)}_\l\vnormf{\xtl}{D}\vnormf{\ytl}{D}-c_m\beta^m\vnormf{\xtlm}{D}\vnormf{\ytl}{D}\\
	&\gtrsim (\alpha^{(3)}_\l-c_m\beta^m)\vnormf{\xtlm}{D}\vnormf{\ytlm}{D}\\
	&\gtrsim \alpha^{(3)}_\l\vnormf{\xtlm}{D}\vnormf{\ytlm}{D}.
	\end{align*}
	In the last inequality, we employed the oversampling condition \eqref{eq:oversampling1}.
	\end{proof}
\end{theorem} 

\begin{remark}[Level oversampling parameter]\label{remark:oversamplingml}
	In the proof of Theorem \ref{lemma:infsuploc}, it is possible to relax oversampling condition \eqref{eq:oversampling1}, choosing level-dependent oversampling parameters	
	\begin{alignat*}{1}
	m_\l\gtrsim|\log(\beta)|^{-1}\cdot\begin{cases}
	|\log(\kappa^{n+1}c_m)|\quad &\text{ for }\l = 1,\\
	|\log(c_m)|\quad &\text{ for }\l\geq 2.
	\end{cases} 
	\end{alignat*}
\end{remark}

\begin{theorem}[A-priori error estimate]\label{theorem:convloc}
	{\color{black}
		Let $u\in\V$ be the solution of \eqref{eq:wfcont} for the right-hand side $f$. If  Assumption \ref{assumption:rescond} and \ref{assumption:poly} are satisfied as well as the oversampling condition \eqref{eq:oversampling1}, then is a constant $c_{\kappa,m,L}$ depending polynomially on $\kappa$, $m$ and linearly on $L$, but independent of $H_L$ such that the solution of the proposed practical multi-level method \eqref{eq:wfdiscmlloc}, \eqref{eq:wfdiscmladduploc} satisfies for any $f\in H^{s}(D)$, $s\in [0,1]$}
	\begin{equation}\label{eq:estconv}
	\vnorm{u-\tilde{u}^m}{D} \lesssim H_L^{1+s}\|f\|_{H^s(D)}+c_{\kappa,m,L}\beta^m\tnorm{f}{D}.
	\end{equation}
	\begin{proof}
		We first estimate using the triangle inequality
		\begin{equation*}
		\vnorm{u-\tilde{u}^m}{D} \leq \vnorm{u-\tilde{u}}{D}+\vnorm{\tilde{u}^m-\tilde{u}}{D}\leq \vnorm{u-\tilde{u}}{D}+\sum_{\l=1}^{L}\vnorm{\xtlm-\xtl}{D}
		\end{equation*}
		with $\xtlm\in\Xtlm$ solving \eqref{eq:wfdiscmlloc}. The first term can be estimated with Lemma \ref{lemma:errideal}. The functions $\xtlm$ are  non-conforming, non-consistent Petrov-Galerkin approximations of $\xtl\in\Xtl$ solving \eqref{eq:wfdiscml}. Using Strang's Lemma (c.f. \cite{Knabner2003}) yields
		\begin{align*}
		\alpha^{(4)}_\l\vnorm{\xtlm-\xtl}{D}&\lesssim\inf_{\xtlm\in \Xtlm}\vnorm{\xtlm-\xtl}{D} +\sup_{\ytlm\in\Ytlm} \frac{\big|a(\xtl,\ytlm)-\tspf{f}{\ytlm}{D}\big|}{\vnormf{\ytlm}{D}}
		\end{align*}
		For the first term, we choose $\xtlm := (1-\C_\l^m)\Pl\xtl$. Using Lemma \ref{lemma:infsupXtlYtl}, we get
		\begin{equation*}
		\inf_{\xtlm\in \Xtlm}\vnorm{\xtlm-\xtl}{D}\leq \vnorm{(\Clm-\C_\l)\Pl\xtl}{D}\lesssim c_m\beta^m\vnorm{\xtl}{D}\lesssim c_m\beta^m{\alpha^{(2),-1}_\l}\tnorm{f}{D}.
		\end{equation*}
		Let $\ytl\in\Ytl$ be arbitrary. Then, for the second term, we obtain after some algebraic manipulations
		\begin{align*}
		a(\xtl,\ytlm)&-\tspf{f}{\ytlm}{D}=a(\xtl,\ytlm)-a(\xtlm,\ytlm)=a(\xtl,\ytl)-a(\xtlm,\ytlm) +\\
		&+a(\xtl,\ytlm-\ytl)=\tspf{f}{\ytl-\ytlm}{D}+a(\xtl,\ytlm-\ytl).
		\end{align*}
		The choice $\ytl := (1-\Cl^*)\Pl\ytlm$ yields
		\begin{align*}
		\big|a(\xtl,\ytlm)-\tspf{f}{\ytlm}{D}\big|&\lesssim \kappa_0^{-1}\tnorm{f}{D}\vnormf{\ytlm-\ytl}{D}+\vnorm{\xtl}{D}\vnormf{\ytlm-\ytl}{D}\\
		&\lesssim \big(\kappa_0^{-1}+{\alpha^{(3),-1}_\l}\big)\tnorm{f}{D}\vnormf{(\Clms-\Clm)\Pl\ytlm}{D}\\
		&\lesssim c_m\beta^m \big(\kappa_0^{-1}+{\alpha^{(3),-1}_\l}\big)\tnorm{f}{D}\vnormf{\ytlm}{D}.
		\end{align*}
		After dividing by $\vnormf{\ytlm}{D}$, the desired estimate follows. 
		Combining both estimates, we obtain 
		\begin{align*}
		\vnorm{\xtlm-\xtl}{D}\lesssim c_{\l,m,\kappa}\beta^m\tnorm{f}{D}\quad\text{ with }\quad c_{\l,m,\kappa}:=c_m\alpha^{(4),-1}_\l\big(\kappa_0^{-1}+\alpha^{(3),-1}_\l\big).
		\end{align*}
		Using Lemma \ref{lemma:infsupXtlYtl}, Remark \ref{infsupl1}, and Theorem \ref{lemma:infsuploc}, we can derive the bounds $c_{1,m,\kappa}\lesssim c_m\kappa^{2n+2}$ for $\l = 1$ and $c_{\l,m,\kappa}\lesssim c_m$ for $\l\geq 2$.
Summing over all levels proves the assertion with the constant
\begin{equation*}
c_{\kappa,m,L} := c_m(\kappa^{2n+2}+L).
\end{equation*}
\end{proof}
\end{theorem}
\begin{remark}[Oversampling]
	If the oversampling parameter, in addition to  \eqref{eq:oversampling1}, also satisfies 
	\begin{equation*}
	m\gtrsim |\log(\beta)|^{-1}|\log(H_Lc_{\kappa,m,L}^{-1})|,
	\end{equation*}
	then the overall error $	\vnorm{u-\tilde{u}^m}{D}$ in \eqref{eq:estconv} is of order $1+s$ in $H_L$.
\end{remark}

{\color{black}
\begin{remark}[$L^2$-error estimate]
	In Theorem \ref{theorem:convloc} an error estimate with respect to the norm $\vnormf{\cdot}{D}$ is stated. Defining $e := u-\tilde u$, an estimate in the weaker $L^2$-norm can be proved for the ideal method
	$$\tnormf{e}{D} = \tnormf{(1-\Pi_L)e}{D}\leq \pi^{-1}H_L \tnormf{\nabla e}{D}\leq \pi^{-1}H_L\vnormf{e}{D},$$
	where we used that $e \in \W_L$. This directly implies  convergence of order $2+s$ for the ideal method \eqref{eq:wfdisc} provided that $f \in H^s(D)$. It is straightforward to extend this result to the localized multi-level method, i.e.,
	$$\|u-\tilde u^m\|_{L^2(D)}\lesssim H_L^{2+s}\|f\|_{H^s(D)} + c_{\kappa,m,L}\beta^m\|f\|_{L^2(D)}.$$
\end{remark}}

{\color{black}
	\begin{remark}[Comparison to classical Helmholtz analysis]
		The analysis of standard finite element discretizations for Helmholtz problems is based on an observation by Schatz \cite{Schatz1974}. He observed that, if an indefinite sesquilinear form satisfies G\r{a}rding's inequality, the stability and quasi-optimality of a finite element discretization can be shown under the two assumptions: (i) approximation properties of the finite element test space and (ii) a sufficiently small mesh size. 
		
		In the original paper introducing the LOD for Helmholtz problems, Schatz's argument is employed for the error analysis of the localized method; see  \cite[Theorem 5.5]{Peterseim2017a}. However, due to the sophisticated multi-resolution structure of the discrete problem, Schatz's argument could not be employed in the present work. Nevertheless, using Strang's lemma and the special properties of the multi-resolution ansatz spaces, we were able to prove similar error estimates also for the multi-resolution method with a slightly larger constant; see Theorem \ref{theorem:convloc}.
		
	\end{remark}
}
\subsection{Practical aspects}
The localized element-corrector problems \eqref{eq:elemcorrdef} are still infinite-dimensional problems and need to be discretized itself. {\color{black} For the ease
of presentation, we restrict ourselves in this paper to the classical case of $\mathcal Q_1$ finite elements. We emphasize that a large variety of discretization schemes can be applied, in particular, also $hp$ adaptive methods.

The patch-local corrector problems \eqref{eq:locelemcorr} are discretized and solved on a fine mesh of the respective patch with mesh-size $h$ satisfying that $\kappa^2h$ is sufficiently small. Under additional geometric requirements, this guarantees stability and quasi-uniformity of the $\mathcal Q_1$ finite element discretization; see \cite{Melenk1995}. }Note that all element-corrector problems are independent and thus, can be solved in parallel. The number of fine-scale problems to be solved on level $\l$ is independent of $\H_\l$ and only depends on the oversampling parameter $m$.
In Lemma \ref{lemma:coercivity}, it was shown that the corrector problems are coercive. However, this is not a real practical benefit, since the spaces $\W_\l$ are difficult to construct numerically. It is more practicable to derive a saddle point formulation of the problem with constraints enforcing that the solution is in $\W_\l$. The number of constraints ($\#$c) is small and only depends polynomially on $m$. 

The solution procedure proposed in \cite{Malqvist2020a} computes $\#$c Helmholtz problems on every (reference) patch. Since the patches have a diameter at most of order $mH_1$, the effective wave number of the patch problems is at most of order $m$ (Assumption \ref{assumption:rescond}). For such Helmholtz problems, there exist effective preconditioners; see \cite{Gander2015,Gander2019,Ramos2020}. With the solution of the Helmholtz patch problems at hand, it is then possible to calculate the Schur complement explicitly. 
For a more detailed discussion of the numerical solution of the corrector problems, see  \cite{Engwer2019,Malqvist2020a}.

\section{Fast solvers}\label{sec:fastsolvers}
In this section, we propose a strategy for solving the level problems \eqref{eq:wfdiscmlloc} efficiently. As Remark \ref{remark:oversamplingml} already indicates, for $\l\geq 2$, problems \eqref{eq:wfdiscmlloc} are good-natured in the sense that it can be solved up to a given tolerance within a fixed number of GMRES steps (Theorem \ref{theorem:uniformsteps}). However, for $\l = 1$, the $\kappa$-dependence in Remark \ref{remark:oversamplingml} makes this impossible. Since this block is relatively small, we suggest a direct solver for the first block. 

Recalling the choice of bases \eqref{eq:basesXtlmYtlm},
we can rewrite \eqref{eq:wfdiscmlloc} for $\l = 1,...,L$ into (decoupled) linear systems of equations
\begin{equation}\label{eq:lse}
\mathbf A^m_\l \mathbf x_\l^m = \mathbf f_\l^m
\end{equation}
with $\mathbf A^m_\l := \big(a(\blkm,\blims)\big)_{j,k=1,...,N_\l}$ and $\mathbf f^m_\l := \big((f,\blims)_{L^2(D)}\big)_{j = 1,...,N_\l}$. 
\subsection{GMRES method}
Since the linear systems of equations \eqref{eq:lse} are sparse but non-hermitian, our preferred iterative solver is the GMRES method; see \cite{Saad2003}.
For general invertible matrices $\mathbf A \in \mathbb{C}^{N\times N}$, its convergence can be expressed in terms of the field of values $\mathcal{F}(\mathbf A)$ which is defined as
	\begin{equation*}
	\mathcal{F}(\mathbf A) := \left\{\frac{\left(\mathbf A \xi,  \xi\right)_{\mathbb{C}^N}}{\left(\xi, \xi\right)_{\mathbb{C}^N}},\; \xi\in \mathbb{C}^N\right\}
	\end{equation*}
	with $\left(\cdot,\cdot\right)_{\mathbb{C}^N}$ denoting the standard inner product on $\mathbb{C}^N$. Henceforth, let $|\cdot|_2$ denote the norm induced by $\left(\cdot,\cdot\right)_{\mathbb{C}^N}$.

Consider a right-hand side $\mathbf f\in \mathbb{C}^N$ and let $\mathbf{x}^{(k)}$ denote the $k$-th iterate of the GMRES method applied to $\mathbf{Ax} = \mathbf f$. Then, the residuals $\mathbf{r}^{(k)}:=\mathbf{A}\mathbf{x}^{(k)}-\mathbf{f}$ converge to zero with the rate  
\begin{equation}\label{eq:gmresconv}
\frac{|\mathbf r^{(k)}|_2}{|\mathbf r^{(0)}|_2}\leq \left(1-\frac{\min_{\mathbf \xi\in \mathcal{F}(\mathbf A)} |\mathbf \xi|_2^2}{\|\mathbf A\|_2^2}\right)^{k/2}.
\end{equation}

This result can be found in \cite{Liesen2020}.


\subsection{Uniform number of GMRES iterations for $\l\geq 2$}
The following theorem relies on uniform upper and lower bounds of the fields of values. With these bounds, the convergence rate of the GMRES method can be uniformly bounded away from one.

\begin{theorem}[Uniform convergence of GMRES applied to \eqref{eq:lse} for $\l\geq 2$]\label{theorem:uniformsteps}
	Suppose that $\l\geq 2$ and  let Assumption \ref{assumption:rescond} be satisfied as well as the oversampling condition 
	\begin{equation}\label{eq:oversampling2}
	m\gtrsim |\log(\beta)|^{-1}|\log(c_m)|.
	\end{equation}
	Then the linear systems \eqref{eq:lse} can be solved up to a specified relative tolerance within a fixed number of iterations depending only on the tolerance.
	\begin{proof}
		Let $\l \in \{2,...,L\}$ be fixed. First, we consider the ideal (non-localized) case, with basis functions $\bli := (1-\Cl)\Pit_\l\phi_{\l,j}$ and $\blis := (1-\Cl^*)\Pit_\l\phi_{\ell,j}$. We derive upper and lower bounds for the field of values of  $\mathbf A_\l := \big(a(\blk,\blis)\big)_{j,k=1,...,N_\l}$. 
		For any $\xi \in \mathbb{C}^{N_\l}$, we obtain the upper bound
		\begin{align*}
		\left(\mathbf A_\l \xi,\xi \right)_{\mathbb{C}^{N_\l}} &= a\Bigg(\sum_{j = 1}^{N_\l}\xi_j \bli,\sum_{j = 1}^{N_l}\xi_j \blis\Bigg) = a\Bigg(\sum_{j = 1}^{N_\l}\xi_j \bli,\sum_{j = 1}^{N_\l}\xi_j \bli\Bigg)\\
		&\lesssim \tnormb{\nabla \Bigg(\sum_{j = 1}^{N_\l} \xi_j\bli\Bigg)}{D}^2= \tnormb{\nabla  (1-\Cl)\Bigg(\sum_{j = 1}^{N_\l} \xi_j \Pit_\l\phi_{\l,j}\Bigg)}{D}^2\\
		&\lesssim \vnormb{\sum_{j = 1}^{N_\l} \xi_j  \Pit_\l\phi_{\l,j}}{D}^2= \kappa^2\tnormb{\sum_{j = 1}^{N_\l} \xi_j \Pit_\l\phi_{\l,j}}{D}^2 + \tnormb{\nabla \Bigg(\sum_{j = 1}^{N_\l} \xi_j  \Pit_\l\phi_{\l,j}\Bigg)}{D}^2\\
		&\lesssim \frac{(H_\l\kappa)^2+1}{H_\l^2}\tnormb{\sum_{j = 1}^{N_\l} \xi_j \phi_{\ell,j}}{D}^2\lesssim \frac{1}{H_\l^2}|\xi|_2^2,
		\end{align*}
		where we used \eqref{eq:Pit}, \eqref{eq:orthogonalitykernel}, as well as Assumption \ref{assumption:rescond}, Lemma \ref{lemma:coercivity} and Remark \ref{remark:wellposednesscorr}. The lower bound can be derived as follows
		\begin{align*}
		\left(\mathbf A_\l \xi,\xi \right)_{\mathbb{C}^{N_\l}} &=a\Bigg(\sum_{j = 1}^{N_\l}\xi_j \bli,\sum_{j = 1}^{N_l}\xi_j \blis\Bigg) = a\Bigg(\sum_{j = 1}^{N_\l}\xi_j \bli,\sum_{j = 1}^{N_l}\xi_j \bli\Bigg)\\
		&\gtrsim \tnormb{\nabla \Bigg(\sum_{j = 1}^{N_l}\xi_j \bli\Bigg)}{D}^2\gtrsim \frac{1}{H_\l^2}\tnormb{(1-\Pi_{\l-1})\Bigg(\sum_{j = 1}^{N_l}\xi_j \bli\Bigg)}{D}^2\\
		& = \frac{1}{H_\l^2}\tnormb{\sum_{j = 1}^{N_l}\xi_j\bli}{D}^2\gtrsim \frac{1}{H_\l^2}\tnormb{\sum_{j = 1}^{N_l}\xi_j\Pi_\l \bli}{D}^2 \\
		&= \frac{1}{H_\l^2}\tnormb{\sum_{j = 1}^{N_l}\xi_j\phi_{\ell,j}}{D}^2 = \frac{1}{H_\l^2}|\xi|^2_2
		\end{align*}
		using \eqref{eq:Pi} and Lemma \ref{lemma:coercivity}. Next, we consider the localized case. With  \eqref{eq:orthogonalitykernel}, we can deduce the relation 
		\begin{align*}
		\left(\mathbf A_\l^m \xi,\xi\right)_{\mathbb{C}^{N_\l}} &= a\Bigg(\sum_{j = 1}^{N_\l}\xi_j \blim,\sum_{j = 1}^{N_l}\xi_j \blims\Bigg) \\
		&= a\Bigg(\sum_{j = 1}^{N_\l}\xi_j \bli,\sum_{j = 1}^{N_l}\xi_j \bli\Bigg)+a\Bigg(\sum_{j = 1}^{N_\l}\xi_j (\blim-\bli),\sum_{j = 1}^{N_l}\xi_j \blims\Bigg).
		\end{align*}
		For the first term, we can use the bounds from the ideal case. The second term is exponentially small, since
		\begin{align*}
		|M_2| &\lesssim  \vnormb{(\Clm-\Cl)\Bigg(\sum_{j = 1}^{N_\l}\xi_j\Pl \Pit_\l\phi_{\l,j}\Bigg)}{D}\vnormb{\left(1-\Clms\right)\Bigg(\sum_{j=1}^{N_\l}\xi_j\Pl \Pit_\l\phi_{\l,j}\Bigg)}{D}\\
		&\lesssim c_m\beta^m\vnormb{\Pl\Bigg(\sum_{j = 1}^{N_\l}\xi_j \Pit_\l\phi_{\l,j}\Bigg)}{D}^2\lesssim \beta^m\vnormb{\sum_{j = 1}^{N_\l}\xi_j \Pit_\l\phi_{\l,j}}{D}^2\\
		&=c_m\beta^m\Bigg(\kappa^2\tnormb{\sum_{j = 1}^{N_\l}\xi_j \Pit_\l\phi_{\l,j}}{D}^2 + \tnormb{\nabla\Bigg(\sum_{j = 1}^{N_\l}\xi_j \Pit_\l\phi_{\l,j}\Bigg)}{D}^2\Bigg)\\
		&\lesssim c_m\beta^m\frac{\left((H_\l\kappa)^2 + 1\right)}{H_\l^2}\tnormb{\sum_{j = 1}^{N_\l}\xi_j \phi_{\l,j}}{D}^2\lesssim \frac{c_m\beta^m}{H_\l^2}|\xi|^2_2,
		\end{align*}
		where we used \eqref{eq:Pit}, \eqref{eq:P}, Assumption \ref{assumption:rescond} and Lemma \ref{lemma:locerrcorr}. The constant  $c_m$ depends polynomially on $m$. Combining the previous estimates yields
		\begin{equation*}
		\frac{\min_{\mathbf \xi\in \mathcal{F}(\mathbf A_\l^m)} |\mathbf \xi|_2^2}{\max_{\mathbf \xi\in \mathcal{F}(\mathbf A_\l^m)} |\mathbf \xi|_2^2}\gtrsim\frac{1-c_m\beta^m}{H_\l^2}\frac{H_\l^2}{1+c_m\beta^m}|\xi|_2^2 = \frac{1-c_m\beta^m}{1+c_m\beta^m}|\xi|_2^2.
		\end{equation*}
		Using oversampling condition \eqref{eq:oversampling2}, we estimate 
		\begin{equation*}
		0<c<\frac{1}{4}\frac{\min_{\mathbf \xi\in \mathcal{F}(\mathbf A_\l^m)} |\mathbf \xi|_2^2}{\max_{\mathbf \xi\in \mathcal{F}(\mathbf A_\l^m)} |\mathbf \xi|_2^2}\leq \frac{\min_{\mathbf \xi\in \mathcal{F}(\mathbf A_\l^m)} |\mathbf \xi|_2^2}{\|\mathbf A_\l^m\|_2^2}<1
		\end{equation*}
		with $c$ independent of $\kappa$, $\l$, $H_\l$, $L$, and $m$. Here, we applied the well-known inequality $$\max_{\mathbf \xi \in \mathcal{F}(\mathbf A_\l^m)}|\mathbf \xi|_2\leq\|\mathbf A_\l^m\|_2 \leq 2 \max_{\mathbf \xi \in \mathcal{F}(\mathbf A_\l^m)}|\mathbf \xi|_2,$$ see \cite[Chapter 1]{Horn1990}. Equation \eqref{eq:gmresconv} yields the assertion.
	\end{proof}
\end{theorem}
\section{Numerical experiments}\label{sec:numexp}
In this section, we present a sequence of numerical results. We begin with a convergence test, which confirms the convergence result from Theorem \ref{theorem:convloc} numerically. 
Next, we demonstrate the improved stability properties (even for the elliptic case)  by the novel construction \eqref{eq:deflocttsp}; see  Remark \ref{remark:stabH1}. 
Moreover, we demonstrate that the proposed method is also applicable to heterogeneous Helmholtz problems. {\color{black} A high-frequency scenario underlines the effectiveness of the proposed numerical method also for such regimes. }
Lastly, we consider a scattering scenario with relatively large wave number and evaluate the condition numbers of the respective levels. This confirms the statement of Theorem \ref{theorem:uniformsteps} numerically.
\subsection{Convergence test}
We consider $D := (0,1)^2$ with $\Gamma_\mathrm{R}:=\partial D$ and the smooth right-hand side 
\begin{equation}\label{eq:sourcesincos}
f(x_1,x_2) = \sin(\pi x_1)\cos(\pi x_2).
\end{equation}
The convergence test is performed for $\kappa = 2^j$, $j = 0,...,3$. For all $\kappa$, we fix a coarse (Cartesian) mesh with $H_1= 2^{-j-1}$ and  perform the numerical computations for mesh hierarchies with $H_L = 2^{-j-1},...,2^{-7}$. The corrector problems \eqref{eq:locelemcorr} are discretized using a fine mesh with mesh size $h = 2^{-9}$. Errors are computed to the reference solution on the same fine mesh. 

\begin{figure}[h]
	\includegraphics[width = \textwidth]{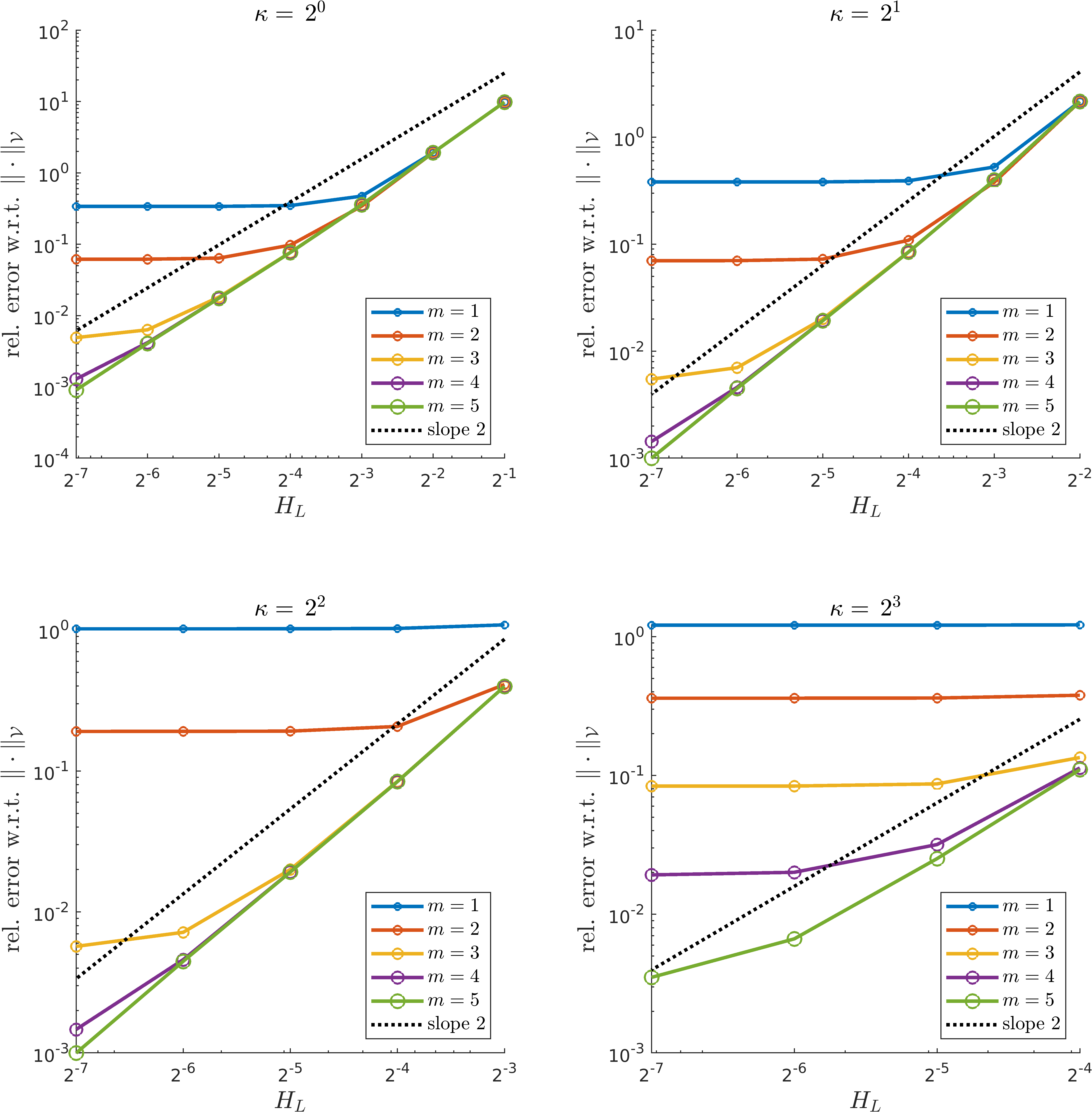}
	\caption{Convergence plots of the proposed multi-level LOD method for oversampling parameters $m = 1,...,5$. From left ro right, top to bottom: $\kappa = 2^j, j= 0,...,3$. For comparison, a line with slope 2 is indicated.}
	\label{fig:conv}
\end{figure}

In the convergence plots of Figure \ref{fig:conv}, one can clearly see the LOD-typical behavior \cite{Malqvist2020a} that the error first decreases with order $2$ ($f\in H^1(D)$), until the localization error dominates, then, the overall error stagnates. This is the expected outcome and is in line with Theorem \ref{theorem:convloc}. For increasing $\kappa$, one observes that for small oversampling parameters $m$, the error does not decrease, but stays of order $\mathcal{O}(1)$. In this case, the oversampling parameter does not fulfill condition \eqref{eq:oversampling1} and convergence cannot be expected. 

\subsection{Improved stability by novel basis construction}\label{subsec:improstab}
A already mentioned in Remark \ref{remark:stabH1}, the novel basis construction \eqref{eq:deflocttsp} (henceforth referred to as \emph{stabilized gamblets}) has superior stability properties compared to previously known constructions (henceforth referred to as \emph{normal gamblets}); see \cite{Owhadi2017a,Owhadi2019,Maier2020a}. Since this is true even for the elliptic case, we use the Poisson problem for demonstration purposes
\begin{alignat*}{1}
-\Delta u &= f\qquad \text{ in }D,\\
u &= 0\qquad \text{ on }\partial D.
\end{alignat*} 
For the numerical experiment, we set $D:=(0,1)^2$ and use the source \eqref{eq:sourcesincos}. For the sake of simplicity, we consider a single-level method ($L=1$) for demonstrating the improved stability properties. For the multi-level method, the observed phenomena are qualitatively the same. We compute numerical approximations for $H_1 = H_L = 2^{-3},...,2^{-7}$ and calculate correctors on a fine mesh with mesh size $h = 2^{-9}$. Errors are calculated to the reference solution on the same fine mesh.

\begin{figure}[h]
	\ContinuedFloat*
	\vspace{.5ex}
		\includegraphics[width=\linewidth]{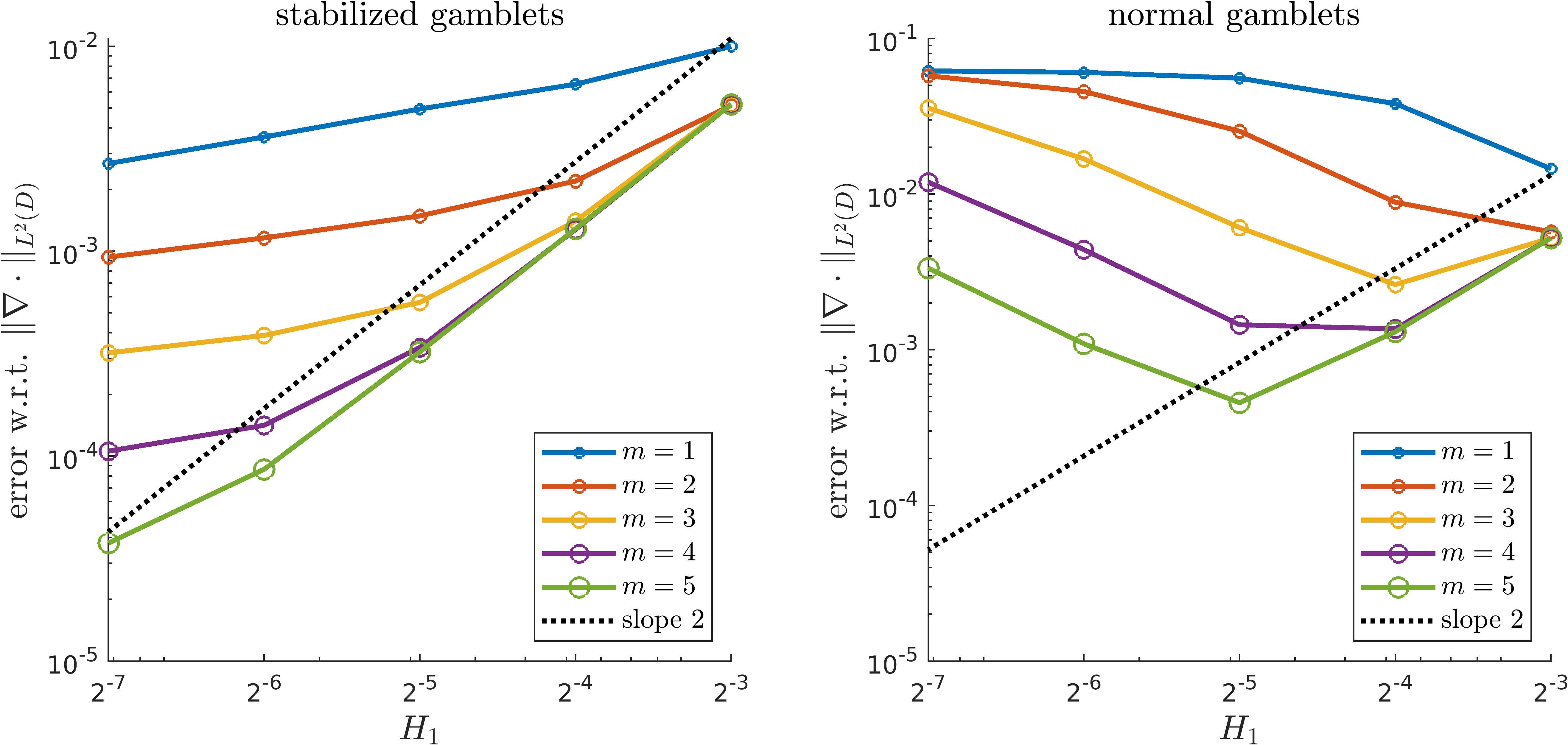}	
		\caption{\label{first} Convergence plots of the single-level method for the Poisson problem for oversampling parameters $m = 1,...,5$. Left: Stabilized gamblets with improved stability properties; Right: Normal gamblets. For comparison, a line with slope 2 is indicated.}
		\label{fig:basisconstrpoisson}
\end{figure}

In Figure \ref{fig:basisconstrpoisson}, one clearly observes the improved stability properties of stabilized gamblets \eqref{eq:deflocttsp} compared to normal gamblets. The error for normal gamblets increases as $H_1$ decreases. Thus, in order to preserve stability, the oversampling parameter $m$ has to be increased as the mesh is refined. However, this is not true for the stabilized gamblets, where we have stability for a fixed oversampling parameter. The error even decreases further with a smaller order. This can be explained by the absence of the block-diagonalization error (c.f. Remark \ref{remark:notbd}) in the single-level case.

For comparison, we also performed this numerical experiment for the Helmholtz problem with $\kappa = 2^3$ and Robin-type boundary conditions. Here, the exact same phenomena can be observed; see Figure \ref{fig:basisconstr}.

\begin{figure}[h]
	\ContinuedFloat
	\includegraphics[width=\linewidth]{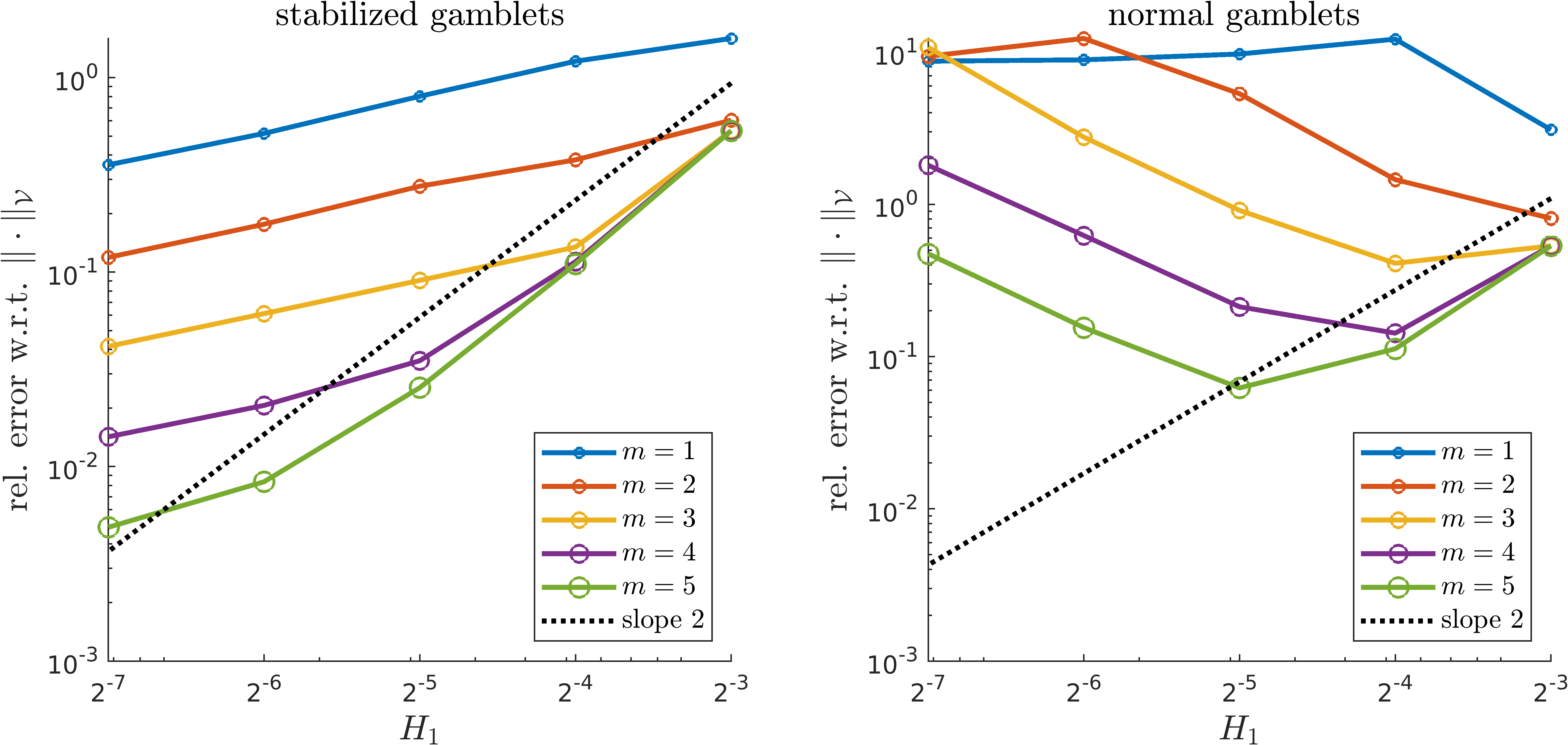}	
	\caption{\label{second} Convergence plots of the single-level method for the Helmholtz problem for oversampling parameters $m = 1,...,5$. Left: Stabilized gamblets with improved stability properties; Right: Normal gamblets. For comparison, a line with slope 2 is indicated.}
	\label{fig:basisconstr}
\end{figure}

\subsection{Variable coefficient}\label{subsec:varcoeff}
Let us consider the Helmholtz problem in heterogeneous media
\begin{equation}\label{eq:varcoeff}
-\nabla\cdot(A\nabla u)-\kappa^2u = f
\end{equation}
with a scalar coefficient $A\in L^{\infty}(D)$ satisfying $0<\gamma\leq A(x)\leq \gamma^\prime <\infty$ for almost all $x\in D$.

\begin{figure}[h]
	\includegraphics[width = \textwidth]{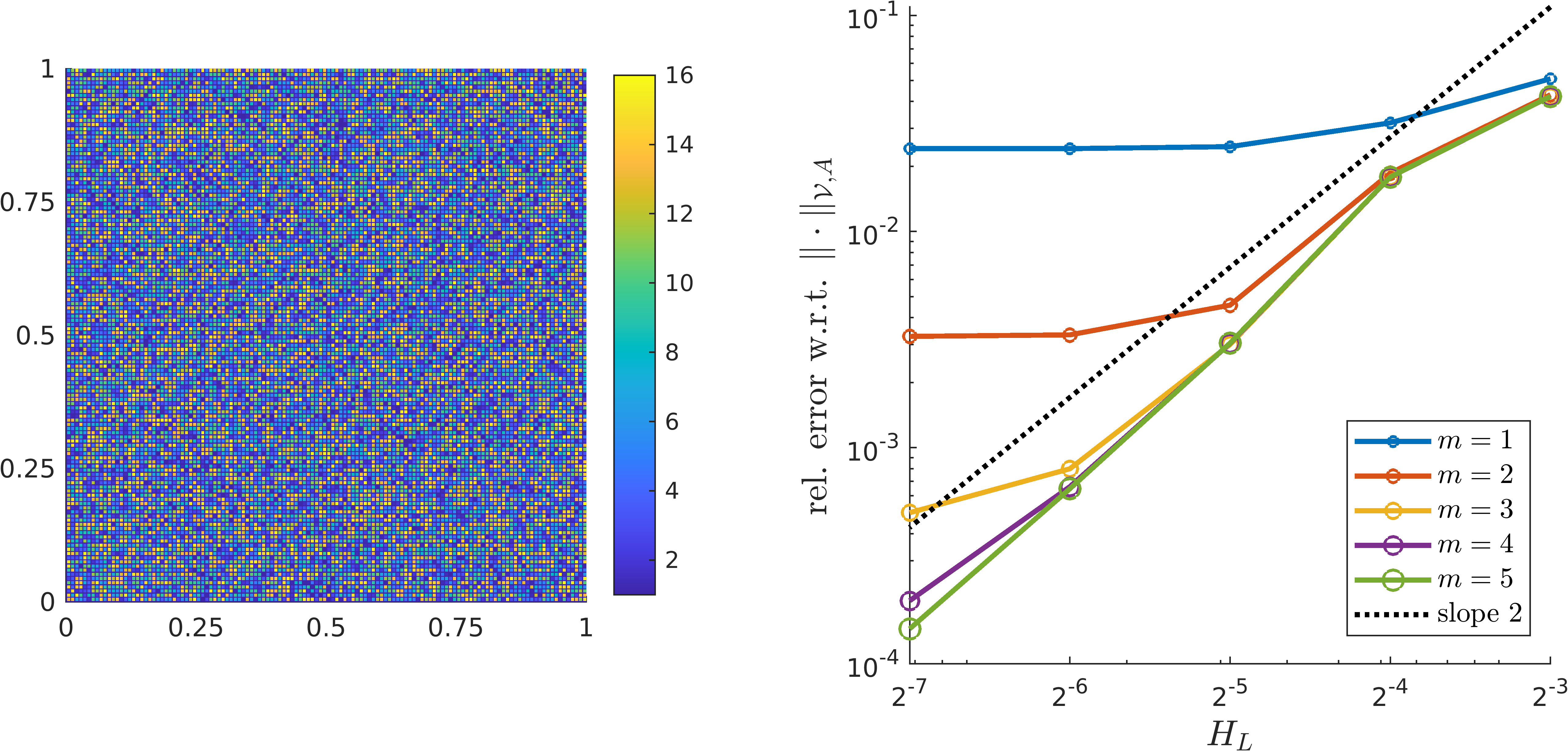}
	\caption{Left: Visualization of the realization of $A$ used for the numerical experiment. Right: Convergence plot of the LOD method for oversampling parameters $m = 1,...,5$. For comparison, a line with slope 2 is indicated.}
	\label{fig:iid}
\end{figure}

For the numerical experiment, we choose $D := (0,1)^2$ with $\Gamma_\mathrm{R}:=\partial D$ and define for $j\in \mathbb{Z}^2$ and the parameter $\epsilon = 2^{-7}$ the inclusions 
\begin{equation*}
S_\epsilon^{j} := \epsilon \left(j+(0.25,0.75)^2\right).
\end{equation*}
On each inclusion in $D$, the coefficient $A$ is chosen to be constant with value uniformly distributed in $[1,16]$. Everywhere else, $A$ is set to $1$. For the numerical experiment, we used the realization of $A$ shown in Figure \ref{fig:iid}. Since $A\vert_{\partial D} = 1$, this choice is compatible with the Robin boundary conditions from \eqref{eq:Helmholtz}. The right-hand side is chosen as 
\begin{alignat}{1}\label{eq:source}
f(x) = \begin{cases}
10000\,\exp\Big(\frac{1}{1-\left(\frac{|x-x_0|}{r}\right)^2}\Big)&\quad \text{ if }\;|x-x_0|<r,\\
0&\quad \text{ else }
\end{cases}
\end{alignat}
with $r = 0.125$  and $x_0 = (0.5,0.5)^T$. In this example, we consider $\kappa = 2^3$. Note that the effective wave number in the inclusions is actually smaller than $\kappa$. As fixed coarse mesh, we use a Cartesian mesh with $H_1 = 2^{-3}$. The numerical computations are performed for hierarchy of meshes with $H_L = 2^{-3},..., 2^{-7}$. The correctors are computed on a fine mesh with $h = 2^{-9}$. Errors are computed to the reference solution on the same fine mesh. As norm for the errors, we choose 
\begin{equation*}
\|u\|_{\V,A}:=\sqrt{\tnorm{A^{1/2}\nabla u}{D}^2+\kappa^2\tnorm{u}{D}^2}.
\end{equation*}

In Figure \ref{fig:iid}, one observes a convergence behavior similar as in Figure \ref{fig:conv}. In contrast to the convergence test for $\kappa= 2^3$ in Figure \ref{fig:conv}, the oversampling parameter $m = 2$ already yields a good approximation in this example. This can be explained by the smaller effective wave number, caused by the coefficient $A$.
{\color{black}
\subsection{High-frequency}
In this numerical experiment, we again choose $D := (0,1)^2$ with pure Robin boundary conditions, i.e.,  $\Gamma_\mathrm{R}:=\partial D$. We consider the wave number $\kappa = 2^6$ and use the right-hand side \eqref{eq:source} with the same parameters as in the previous experiment. As fixed coarse mesh, we use a Cartesian mesh with $H_1 = 2^{-6}$. The numerical computations are performed for hierarchy of meshes with $H_L = 2^{-6},..., 2^{-8}$. The correctors are computed on a fine mesh with $h = 2^{-10}$. Errors are computed to the reference solution on the same fine mesh. 
\begin{figure}[h]
	\includegraphics[width=.5\linewidth]{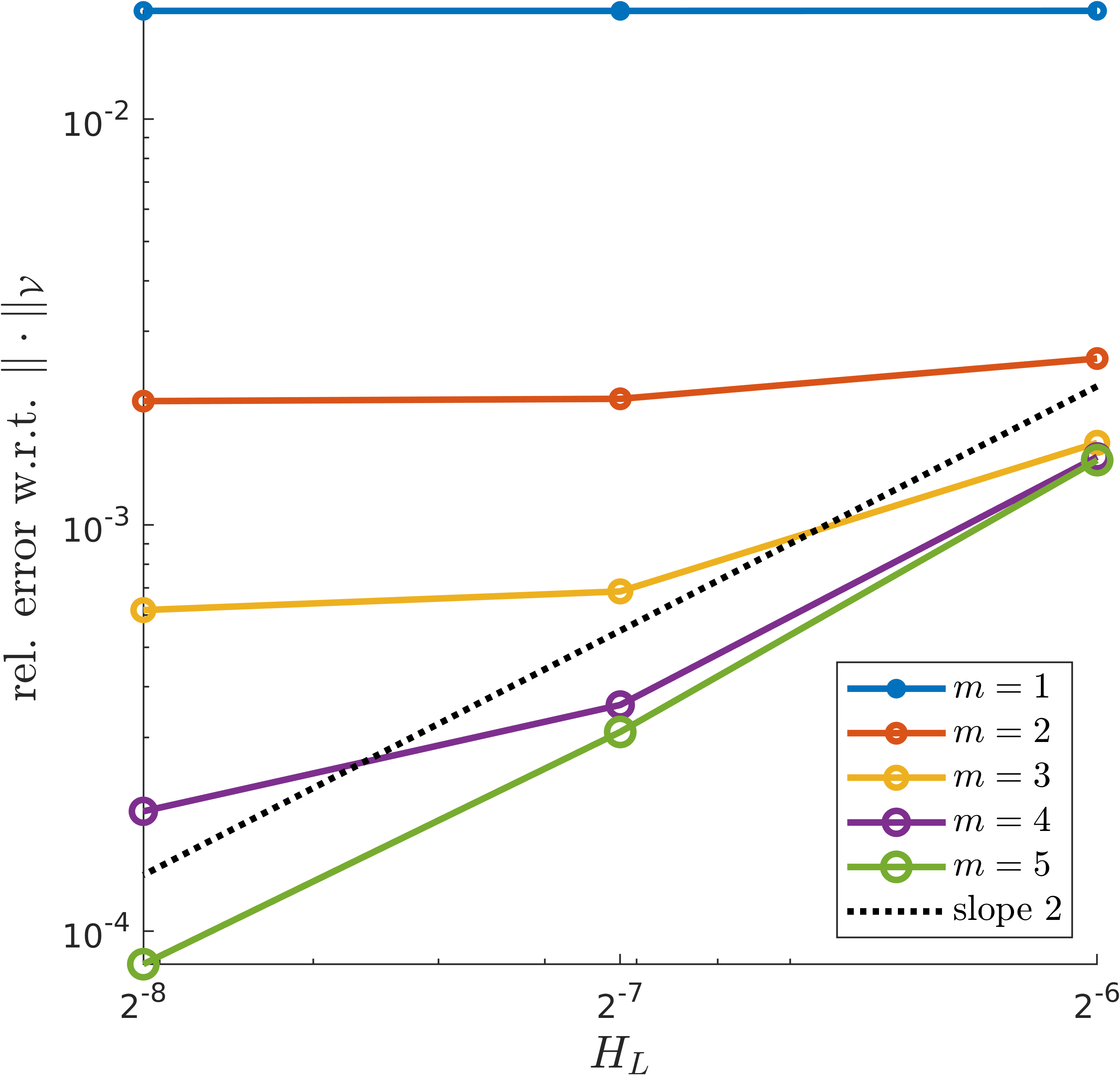}
	\caption{ Convergence plots of the proposed multi-level LOD method for the high-frequency example with $\kappa = 2^6$ for oversampling parameters $m = 1,...,5$. For comparison, a line with slope 2 is indicated.}
	\label{fig:hf}
\end{figure}

In Figure \ref{fig:hf} one can observe the convergence behavior as expected from Theorem \ref{theorem:convloc}, i.e., if the oversampling parameter $m$ is chosen sufficiently large (i.e., $m \gtrsim \log(\kappa)$), one has $\kappa$-independent convergence of the proposed multi-resolution method. This numerical example proves the effectiveness of the numerical method also for the high-frequency regime.
}
\subsection{Uniform number of GMRES steps}
\begin{figure}[h]
	\includegraphics[width = .6\textwidth]{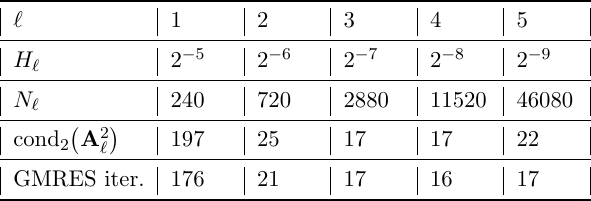}
	\caption{
		Table with properties of the linear system of equations for levels $\l = 1,...,5$.}
	\label{fig:scattering}
\end{figure}

We consider $D := (0,1)^2\backslash S$ with the scatterer $S=[0.375,0.625]^2$. At the boundary of $S$, we impose Dirichlet boundary conditions $\Gamma_\mathrm{D}=\partial S$ and  Robin boundary conditions are imposed at the artificial boundary $\Gamma_\mathrm{R}:=\partial D\backslash\partial S$. We use the source \eqref{eq:source} from above with $r = 0.05$ and  $x_0 = (0.125,0.125)$. We choose $\kappa = 2^5$. The hierarchy used for the computations is specified in Figure \ref{fig:scattering}. The correctors are computed on a fine mesh with $h = 2^{-10}$. The oversampling parameter is chosen as $m = 2$. We use the GMRES method with restart after $50$ iterations. The GMRES iteration terminates if a relative residual of $10^{-6}$ is reached.

The table in Figure \ref{fig:scattering} clearly shows that uniform boundedness of the number of GMRES iterations for $\l \geq 1$. As expected, the GMRES iterations performs very poorly for $\l = 1$, which suggests to use a direct solver for the relatively small linear system of equations corresponding to $\l = 1$. The respective condition numbers, which are a good indicator for how fast iterative solvers converge, underline this observation.

%

\bibliographystyle{alpha}
\bibliography{bib}

\appendix
\section{}\label{sec:appendix}
\begin{proof}[Proof of Lemma \ref{lemma:infsupUtVt}]
			This proof is based on the proof of \cite[Theorem 4.4]{Peterseim2017a}. Note that $(1-\C_L^*)$ is a Fortin operator (as in the theory of mixed methods \cite{Fortin1977}), since for all $\tilde{u}\in \Ut$, $v\in\V$
			\begin{equation*}
			a(\tilde{u},(1-\C_L^*)v) = a(\tilde{u},v) - a(\tilde{u},\C_L^*v) = a(\tilde{u},v),
			\end{equation*}
			where we used \eqref{eq:defcorr}. The continuity of $(1-\C_L^*)$ follows with Remark \ref{remark:wellposednesscorr}. 
			\begin{align*}
			&\inf_{\tilde u\in \Ut}\sup_{\tilde v\in\Vt} \frac{\mathfrak{R}a(\tilde{u},\tilde{v})}{\vnorm{\tilde{u}}{D}\vnorm{\tilde{v}}{D}} = \inf_{\tilde u\in \Ut}\sup_{v\in\V} \frac{\mathfrak{R}a(\tilde{u},(1-\C_L^*)v)}{\vnorm{\tilde{u}}{D}\vnorm{(1-\C_L^*)v}{D}}\\
			& \gtrsim \inf_{\tilde u\in \Ut}\sup_{v\in\V} \frac{\mathfrak{R}a(\tilde{u},v)}{\vnorm{\tilde{u}}{D}\vnorm{v}{D}}
			\geq \inf_{u\in \U}\sup_{v\in\V} \frac{\mathfrak{R}a({u},v)}{\vnorm{{u}}{D}\vnorm{v}{D}}\gtrsim \alpha^{(1)}.
			\end{align*}
			Thus $\alpha^{(2)}$ has the same $\kappa$-dependence as $\alpha^{(1)}$, i.e., $\alpha^{(2)} \simeq \alpha^{(1)}\simeq ( c_{\mathrm{stab}}(\kappa)\kappa)^{-1}$.
\end{proof}
\begin{proof}[Proof of Lemma \ref{lemma:errideal}]
	For  $s = 0$, the proof can also be found in \cite[Theorem 4.5]{Peterseim2017a}.
	Due to Galerkin orthogonality, we have that $e := u - \tilde u = \C_L u \in \W_L$. This implies
	\begin{align*}
	\frac{1}{2}\vnorm{e}{D}^2 &\leq \mathfrak{R}a(e,e) = \mathfrak{R}\tsp{f}{e}{D} =\mathfrak{R}\tsp{f-\Pi_Lf}{e}{D}\\
	&\leq \tnorm{f-\Pi_Lf}{D}\tnorm{(1-\Pi_L)e}{D}\leq \pi^{-1}  H_L \tnorm{f-\Pi_Lf}{D}\vnorm{e}{D}.
	\end{align*}
	If $f\in H^1(D)$, then \eqref{eq:Pi} yields an additional order in $H_L$.
	Dividing by $\vnorm{e}{D}$, the result follows for $s\in\{0,1\}$. By arguments from interpolation theory (see e.g. \cite {Brenner2008}), the assertion can be concluded for any $s\in [0,1]$ {\color{black}with a non-explicit constant $c$.} 
\end{proof}

\begin{proof}[Proof of Lemma \ref{lemma:expdec}]
	Let $\l\in\{1,...,L\}$ be fixed. For shorter notation denote $\varphi:= \Clt v$. Define the finite element cut-off function $\eta\in W^{1,\infty}(D,[0,1])$ such that
	\begin{alignat*}{2}
	&\eta \equiv 0 \quad &&\text{ in }N_\l^{m-1}(T),\\
	&\eta \equiv 1 \quad &&\text{ in }D\backslash N_\l^{m}(T),\\
	&0 \leq \eta \leq 1 \qquad &&\text{ in } R:=N_\l^{m}(T)\backslash N_\l^{m-1}(T).
	\end{alignat*}
	{\color{black} 
	with $\|\nabla \eta\|_{L^\infty(D)}\leq H_\l^{-1}$. We estimate
	\begin{align*}
	\tnorm{\nabla \varphi}{D\backslash N_\l^m(T)}^2 &= \mathfrak{R} \tsp{\nabla\varphi}{\nabla\varphi}{D\backslash N_\l^m(T)}\leq\mathfrak{R} \tsp{\nabla\varphi}{\eta\nabla\varphi}{D}\\
	&= \mathfrak{R}\tsp{\nabla \varphi}{\nabla (\eta\varphi)}{D} - \mathfrak{R}\tsp{\nabla \varphi}{\varphi \nabla \eta}{D}\\
	&\leq \Big\vert \mathfrak{R}\tsp{\nabla \varphi}{\nabla (\eta\varphi - \widetilde\Pi_\l\Pi_\l(\eta\varphi))}{D}\Big\vert   \\
	&\qquad +\Big\vert \mathfrak{R}\tsp{\nabla \varphi}{\nabla\widetilde{\Pi}_\l\Pi_\l(\eta\varphi)}{D}\Big\vert + \big\vert\mathfrak{R}\tsp{\nabla \varphi}{\varphi \nabla \eta}{D}\big\vert\\
	&=: M_1 + M_2 + M_3.
	\end{align*}
	For $M_1$, we have $w:=(1-\widetilde{\Pi}_\l\Pi_\l)(\eta\varphi) \in \W_\l$ and $w\equiv 0$ in $T$ for $m\geq 1$, thus $a_T(\varphi,w) = 0$. Furthermore it holds $\widetilde{\Pi}_\l\Pi_\l(\eta\varphi) \equiv 0$ in $D\backslash N_\l^{m}(T)$ since $\varphi\in\W_\l$. Using the definition of the corrector problems \eqref{eq:elemcorrdef}, \eqref{eq:Pi} and \eqref{eq:Pit}, we get
	\begin{align*}
	M_1 &= \Big\vert \mathfrak{R}\tsp{\nabla \varphi}{\nabla (\eta\varphi - \widetilde\Pi_\l\Pi_\l(\eta\varphi))}{D}\Big\vert\\
	&=\Big\vert\mathfrak{R}\kappa^2  \tsp{\varphi}{\eta\varphi-\widetilde{\Pi}_\l\Pi_\l(\eta\varphi)}{D} \Big\vert\\
	&\leq \kappa^2 \big\vert \tsp{\varphi}{\eta\varphi}{D}\big\vert + \kappa^2\Big\vert\tsp{\varphi}{\widetilde{\Pi}_\l\Pi_\l(\eta\varphi)}{D}\Big\vert\\
	&\leq \kappa^2 \tnorm{(1-\Pi_\l)\varphi}{D\backslash N_\l^{m-1}(T)}^2 + \kappa^2\tnorm{(1-\Pi_\l)\varphi}{R}^2\\
	&\leq \pi^{-2}(\kappa H_\l)^2 \tnorm{\nabla\varphi}{D\backslash N_\l^{m-1}(T)}^2 + \pi^{-2}(\kappa H_\l)^2 \tnorm{\nabla\varphi}{R}^2\\
	&\leq \frac{1}{2}\tnorm{\nabla\varphi}{D\backslash N_\l^{m}(T)}^2 + \frac{1}{2} \tnorm{\nabla\varphi}{R}^2,
	\end{align*}
	where we used Assumption \ref{assumption:rescond}. The left term can be absorbed into $\tnorm{\nabla\varphi}{D\backslash N_\l^m(T)}^2$.
	\\[.5ex]
	Using that $\widetilde{\Pi}_\l\Pi_\l(\eta\varphi) \equiv 0$ in $D\backslash N_\l^{m}(T)$, \eqref{eq:Pi} and \eqref{eq:Pit}, we obtain for $M_2$
	\begin{align*}
	M_2 &= \Big\vert \mathfrak{R}\tsp{\nabla \varphi}{\nabla\widetilde{\Pi}_\l\Pi_\l(\eta\varphi)}{D}\Big\vert\leq c_\mathrm{inv}H_\l^{-1}\tnorm{\nabla\varphi}{R}\tnorm{\eta\varphi}{R}\\
	&\leq c_\mathrm{inv}H_\l^{-1}\tnorm{\nabla\varphi}{R}\tnorm{(1-\Pi_\l)\varphi}{R}\leq c_\mathrm{inv}\pi^{-1} \tnorm{\nabla\varphi}{R}^2.
	\end{align*}
	For $M_3$, we obtain similarly 
	\begin{align*}
	M_3 &=  \big\vert\mathfrak{R}\tsp{\nabla \varphi}{\varphi \nabla \eta}{D}\big\vert \leq H_\l^{-1}\tnorm{\nabla \varphi}{R}\tnorm{(1-\Pi_\l)\varphi}{R}\leq \pi^{-1} \tnorm{\nabla\varphi}{R}^2.
	\end{align*}
}Combining the estimates for $M_1$, $M_2$, and $M_3$, yields
	\begin{align*}
	\frac{1}{2}\tnorm{\nabla\varphi}{D\backslash N^m(T)}^2\leq c \tnorm{\nabla\varphi}{R}^2 &=   c \tnorm{\nabla\varphi}{D\backslash N^{m-1}(T)} - c \tnorm{\nabla\varphi}{D\backslash N^m(T)}^2\\
	\Leftrightarrow\; \Big(\frac{1}{2}+c\Big)\tnorm{\nabla\varphi}{D\backslash N^m(T)}^2&\leq c \tnorm{\nabla\varphi}{D\backslash N^{m-1}(T)}^2\\
	\Leftrightarrow \; \tnorm{\nabla\varphi}{D\backslash N^m(T)} &\leq \beta \tnorm{\nabla\varphi}{D\backslash N^{m-1}(T)},
	\end{align*}
	where $\beta := \sqrt{\frac{c}{\frac{1}{2}+c}}<1$. Iterating this argument, the assertion follows.
\end{proof}
For the proof of Lemma \ref{lemma:locerrcorr}, we infer the intermediate result
\begin{lemma}[Localization error of element-correctors]\label{lemma:locerrelemcorr}
	
	If Assumption \ref{assumption:rescond} is satisfied, then it holds for all $\l\in\{1,...,L\}$, $T\in\T_\l$, $v\in \V$ and $m \in\mathbb{N}$ that
	\begin{equation*}
	\tnorm{\nabla(\Cltm v-\Clt v)}{D} \lesssim \beta^m \tnorm{\nabla \Clt v}{D},
	\end{equation*}
	where $\beta$ is the constant from Lemma \ref{lemma:expdec}. An analogous result holds for $\Cltms-\Clt^*$.
	\begin{proof}
		Let $\l\in\{1,...,L\}$ be fixed. Moreover, let $\varphi:=\Clt v$ and $\varphi^m:=\Cltm v$ and let $w_\l\in \W_\l^m(T)$ be  arbitrary. We obtain, using Lemma \ref{lemma:coercivity} and the Galerkin orthogonality
		\begin{align*}
		\tnorm{\nabla (\varphi^m - \varphi)}{D}^2&\lesssim \mathfrak{R}a(\varphi^m - \varphi,\varphi^m - \varphi) = \mathfrak{R}a(\varphi^m - \varphi,w_\l - \varphi)\\
		&\lesssim\tnorm{\nabla (\varphi^m - \varphi)}{D}\tnorm{\nabla(w_\l - \varphi)}{D}.
		\end{align*}
		Define $w_\l := (1-\widetilde{\Pi}_\l\Pi_\l)(\eta\varphi)$ with 
		the  finite element cut-off function $\eta\in W^{1,\infty}(D,[0,1])$ such that 
		\begin{alignat*}{2}
		&\eta \equiv 1 \quad &&\text{ in }N_\l^{m-1}(T),\\
		&\eta \equiv 0 \quad &&\text{ in }D\backslash N_\l^{m}(T),\\
		&0 \leq \eta \leq 1 \qquad &&\text{ in } R:= N_\l^{m}(T)\backslash N_\l^{m-1}(T)
		\end{alignat*}
		and $\|\nabla \eta\|_{L^\infty(D)}\leq H_\l^{-1}$. Since $w_\l \in \W_\l^m(T)$, we obtain
		\begin{align*}
		\tnorm{\nabla (\varphi^m - \varphi)}{D}&\lesssim \tnorm{\nabla (w_\l - \varphi)}{D} \\
		&\lesssim\tnormf{\nabla ((1-\widetilde{\Pi}_\l\Pi_\l)(1-\eta)\varphi)}{D\backslash N^{m-1}_\l(T)}\\
		&\lesssim\tnorm{\nabla \varphi}{D\backslash N^{m-1}_\l(T)}\\
		&\lesssim\beta^{m-1} \tnorm{\nabla\varphi}{D}\\
		&\lesssim\beta^m\tnorm{\nabla\varphi}{D},
		\end{align*}
		where we used that
		\begin{align*}
		&\tnormf{\nabla ((1-\widetilde{\Pi}_\l\Pi_\l)(1-\eta)\varphi)}{D\backslash N^{m-1}_\l(T)}\\
		&\quad \lesssim \tnorm{\nabla ((1-\eta)\varphi)}{D\backslash N^{m-1}_\l(T)} + H_\l^{-1}\tnorm{ (1-\eta)\varphi}{D\backslash N^{m-1}_\l(T)}\\
		&\quad \lesssim \tnorm{\nabla \varphi}{D\backslash N^{m-1}_\l(T)} + \tnorm{\varphi\nabla\eta}{R} + H_\l^{-1}\tnorm{\varphi}{D\backslash N^{m-1}_\l(T)}\\
		&\quad \lesssim \tnorm{\nabla \varphi}{D\backslash N^{m-1}_\l(T)}.
		\end{align*}
		Here, we applied \eqref{eq:Pi} and \eqref{eq:Pit} multiple times.
	\end{proof}
\end{lemma}

\begin{proof}[Proof of Lemma \ref{lemma:locerrcorr}]
		Let $\l\in\{1,..,L\}$ be fixed. Define $z:=\Clm v - \Cl v$ and $z_T = \Cltm v- \Clt v$ for $T\in\T_\l$. It holds
		\begin{equation*}
		\tnorm{\nabla z}{D}^2\lesssim a(z,z) = \sum_{T\in \mathcal{T}_\ell} a(z_T,z).
		\end{equation*}
		For $T\in\T_\l$ define finite element the cut-off function $\eta\in W^{1,\infty}(D,[0,1])$ such that
		\begin{alignat*}{2}
		&\eta \equiv 0 \quad &&\text{ in }N_\l^{m}(T),\\
		&\eta \equiv 1 \quad &&\text{ in }D\backslash N_\l^{m+1}(T),\\
		&0 \leq \eta \leq 1 \qquad &&\text{ in } R:=N_\l^{m+1}(T)\backslash N_\l^{m}(T)
		\end{alignat*}
		and $\|\nabla \eta\|_{L^\infty(D)}\leq H_\l^{-1}$. Using $\supp((1-\Pit_\l\Pi_\l)(\eta z))\subset D\backslash N_\l^m(T)$,  $(1-\Pit_\l\Pi_\l)(\eta z)\in W_\l$, and \eqref{eq:elemcorrdef}, we get
		\begin{equation*}
		a(z_T,(1-\Pit_\l\Pi_\l)(\eta z)) = -a(\Clt v,(1-\Pit_\l\Pi_\l)(\eta z)) = 0.
		\end{equation*}
		Hence, using that $z = (1-\Pit_\l\Pi_\l )z$, yields
		\begin{equation*}
		a(z_T,z) = a(z_T,z-(1-\Pit_\l\Pi_\l)(\eta z)) = a(z_T,(1-\Pit_\l\Pi_\l )((1-\eta)z)).
		\end{equation*}
		With \eqref{eq:Pi} and \eqref{eq:Pit}, one obtains the bound	
		\begin{equation*}
		\big\vert a(z_T,(1-\Pit_\l\Pi_\l )((1-\eta)z))\big\vert \lesssim \vnorm{z_T}{D}\vnorm{z}{N_\l^{m+1}(T)}.
		\end{equation*}
		The  element-correctors satisfy the  estimate $\tnorm{\nabla \Clt v}{D}\lesssim \vnorm{v}{T}$, since 
		\begin{align*}
		\tnorm{\nabla \Clt v}{D}^2\lesssim a(\Clt v,\Clt v) = a_T(v,\Clt v)
		\lesssim \vnorm{v}{T}\tnorm{\nabla \Clt v}{D}.
		\end{align*}
		Using this and Lemma \ref{lemma:coercivity}, Lemma \ref{lemma:locerrelemcorr} and the finite overlap of the patches, we get after summing over all elements 
		\begin{align*}
		\tnorm{\nabla z}{D}^2&\lesssim \sum_{T\in \mathcal{T}_\ell} a(z_T,z)\lesssim\sum_{T\in \mathcal{T}_\ell} \vnorm{z_T}{D}\vnorm{z}{N_\l^{m+1}(T)}\\
		&\lesssim \beta^m \sum_{T\in \mathcal{T}_\ell}\tnorm{\nabla \Clt v}{D}\vnorm{z}{N_\l^{m+1}(T)}\\
		&\lesssim \beta^m\sum_{T\in \mathcal{T}_\ell} \vnorm{v}{T}\vnorm{z}{N_\l^{m+1}(T)}\\
		&\lesssim c_m\beta^m  \vnorm{v}{D}\vnorm{z}{D}\lesssim c_m\beta^m  \vnorm{v}{D}\tnorm{\nabla z}{D}.
		\end{align*}
		with constant $c_m$ depending polynomially on the oversampling parameter $m$.
\end{proof}
\end{document}